\documentclass[12pt, a4paper]{article}

\usepackage{amsmath,amssymb}

\numberwithin{equation}{section}
\newtheorem{dfn}{Definition}[section]
\newtheorem{thm}[dfn]{Theorem}
\newtheorem{rmk}[dfn]{Remark}
\newtheorem{lem}[dfn]{Lemma}
\newtheorem{cor}[dfn]{Corollary}
\newtheorem{emp}[dfn]{Example}

\newenvironment{prf}{\noindent \textit{Proof} \ }{\hfill $\Box$}

\newcommand{\N}{\mathbb{N}}
\newcommand{\Z}{\mathbb{Z}}
\newcommand{\R}{\mathbb{R}}

\newcommand{\A}{\mathcal{A}}
\newcommand{\F}{\mathcal{F}}
\newcommand{\E}{\mathcal{E}}
\newcommand{\CL}{\mathcal{L}}
\renewcommand{\AA}{\hat{\mathcal{A}}}
\newcommand{\FF}{\hat{\mathcal{F}}}
\newcommand{\EE}{\hat{\mathcal{E}}}
\newcommand{\nn}{\nonumber}

\newcommand{\p}{\partial}
\newcommand{\e}{\epsilon}
\newcommand{\hM}{\hat{M}}
\newcommand{\hU}{\hat{U}}

\DeclareMathOperator{\Ker}{Ker}
\DeclareMathOperator{\Img}{Im}
\DeclareMathOperator{\Der}{Der}

\begin{document}

\title{Bihamiltonian Cohomologies and Integrable Hierarchies I: A Special Case}
\author{Si-Qi Liu, Youjin Zhang\\
{\small Department of Mathematical Sciences,
Tsinghua University}\\
{\small Beijing 100084, P. R. China}}
\date{}\maketitle

\begin{abstract}
We present some general results on properties of the bihamiltonian cohomologies associated to bihamiltonian structures
of hydrodynamic type, and compute the third cohomology for the bihamiltonian structure of the dispersionless KdV hierarchy. 
The result of the computation enables us to prove the existence of  bihamiltonian deformations of the dispersionless KdV
hierarchy starting from any of its infinitesimal deformations. 
\end{abstract}

\tableofcontents

\section{Introduction}

The notion of bihamiltonian cohomologies was introduced in \cite{DZ} for the study of deformations of bihamiltonian structures of hydrodynamic type and the associated integrable hierarchies. Such a class of bihamiltonian integrable hierarchies has important applications in the study of Gromov-Witten theory, singularity theory and some other research fields in mathematical physics, see \cite{Wit, Kon, Dij, Egu, Du1, DZ-cmp, DZ, Gi-2} and references therein. For any semisimple bihamiltonian structure of hydrodynamic type, the knowledge of the first and second  bihamiltonian cohomologies are used in \cite{LZ1, DLZ1} to classify their infinitesimal deformations. The purpose of the present and the subsequent paper \cite{DLZ2} is to study the existence problem of deformations of semisimple bihamiltonian structures of hydrodynamic type with a given infinitesimal deformation. The main tool of our study is again provided by the bihamiltonian cohomologies.

In order to give readers a quick introduction of the notion of bihamiltonian cohomologies, let us first recall the notion of
Poisson cohomology defined by Lichnerowicz for finite dimensional Poisson manifolds \cite{Lich-1}, and then    
introduce the notion of bihamiltonian cohomologies for a Poisson manifold endowed with a second compatible Poisson structure.

Let $M$ be an $n$-dimensional smooth manifold, and let 
\[\Lambda^*=\Lambda^{0}\oplus \Lambda^1\oplus \Lambda^2\oplus\dots\]
be the space of multi-vectors on $M$.
On $\Lambda^*$ there is defined a Schouten-Nijenhuis bracket 
\[[\ \,,\ ]:\Lambda^p\times\Lambda^q \to \Lambda^{p+q-1},\quad (P, Q)\mapsto [P, Q]\]
which makes $(\Lambda^*, [\ \,,\ ])$ a graded Lie algebra (after adjusting the gradation and signs).
Suppose $M$ is endowed with a Poisson bivector $P$, i.e. $P\in\Lambda^2$ such that $[P, P]=0$, then $M$ is called a Poisson manifold and the adjoint
action $d=\mathrm{ad}_P$ of $P$ is a differential on $\Lambda^*$, so $(\Lambda^*, [\ \,,\ ], d)$ becomes a differential
graded Lie algebra (DGLA). It is well-known that this DGLA controls the formal deformation theory of the Poisson manifold
$(M, P)$.

A deformation of $P$ is a formal power series
\begin{equation}\label{ky-2}
\tilde{P}=P+\xi\,Q_1+\xi^2\,Q_2+\cdots
\end{equation}
such that $[\tilde{P}, \tilde{P}]=0$. If we denote the perturbation part $\tilde{P}-P$ by $Q$, then $Q$ satisfies the
Maurer-Cartan equation
\[d\,Q+\frac12[Q,Q]=0\]
which is equivalent to a series of equations 
\[d\,Q_1=0,\ d\,Q_2+\frac12[Q_1,Q_1]=0, \ \dots\]
for $Q_j,\, j\ge 1$.
The Maurer-Cartan equation admits a group of symmetries
\begin{equation}\label{ky-1}
\tilde{P}\mapsto \exp(\xi\, \mathrm{ad}_X)\tilde{P},\quad X\in\Lambda^1,
\end{equation}
which are called gauge transformations. The classification problem of deformations of $P$ is to ask
what is the orbit space of the space of solutions of the Maurer-Cartan equation under the action of gauge transformations.
In particular, deformations in the orbit of $P$ are called trivial deformations.

A solution $Q_1\in\Lambda^2$ of $d Q_1=0$ yields a deformation $P+\xi\, Q_1$ of the Poisson structure $P$ at the approximation up to the first order of $\xi$, and such a bivector $Q_1$ is called an infinitesimal deformation of $P$. From the form of the symmetries \eqref{ky-1} of the Maurer-Cartan equation it follows that an infinitesimal deformation $Q_1$ 
is trivial if and only if there exists a vector field $X\in\Lambda^1$ such that $Q_1=d\,X$. So the infinitesimal
deformations are classified by the elements of the following quotient space:
\[ H^2(M, P)=\frac{\left.\Ker d\right|_{\Lambda^2}}{\left.\Img d\right|_{\Lambda^{1}}}.\]
This space is called the second Poisson cohomology of $(M, P)$. More generally, one can define the $k$-th Poisson cohomology
of $(M, P)$ as follows:
\[ \left. H^0(M,P)=\Ker d\right|_{\Lambda^0},\quad H^k(M, P)=\frac{\left.\Ker d\right|_{\Lambda^k}}{\left.\Img d\right|_{\Lambda^{k-1}}},\quad  k\ge 0.\]
Given an infinitesimal deformation $Q_1$ of the Poisson structure $P$, we need to consider 
whether it is possible to extend it to a deformation \eqref{ky-2} of $P$. The first step is to solve 
the following equation for $Q_2$: 
\[d\,Q_2=-\frac12[Q_1, Q_1]=:W\in\Lambda^3.\]
Note that the right hand side of the above equation satisfies $d\,W=0$, so it defines a cohomology class in $H^3(M, P)$.
If this cohomology group vanishes, then we can solve the above equation to obtain $Q_2$. Further more,
by a similar argument one can show that if $H^3(M, P)$ vanishes then we can always obtain a deformation \eqref{ky-2} of $P$ starting from any of its infinitesimal deformation $Q_1$.
For more details on the relationship between DGLA and formal deformation theory, see e.g. \cite{MANETTI} and references therein.

Now let us assume that the Poisson manifold $(M, P)$ is endowed with a second Poisson bivector $P_2$ which is compatible with $P_1=P$, i.e.  $[P_1, P_2]=0$. In this case we say that the manifold $M$ has a bihamiltonian structure $(P_1, P_2)$.
The two Poisson structures define on $\Lambda^*$ two differentials
\[d_1=\mathrm{ad}_{P_1},\quad  d_2=\mathrm{ad}_{P_2}\]
which satisfy
\[d_1^2=d_1 d_2+d_2 d_1=d_2^2=0.\]
It is easy to see that the kernel of $d_1$, denoted by $\tilde{\Lambda}^*=\Ker d_1$, is a graded Lie subalgebra of
$\Lambda^*$, and $d_2$ is still a differential on $\tilde{\Lambda}^*$, so we have another
differential graded Lie algebra
$(\tilde{\Lambda}^*, [\ \,,\ ], d_2)$.

Similar to the definition of deformation \eqref{ky-2} of a Poisson structure, we can also consider deformations
of a bihamiltonian structure $(P_1, P_2)$ of the form
\begin{equation}
 (\tilde{P}_1, \tilde{P}_2)=(P_1+\xi Q_{1,1}+\xi^2 Q_{1,2}+\dots, P_2+\xi Q_{2,1}+\xi^2 Q_{2,2}+\dots),
\end{equation}
which are required to satisfy the following equations:
\[[\tilde{P}_1, \tilde{P}_1]=[\tilde{P}_1, \tilde{P}_2]=[\tilde{P}_2, \tilde{P}_2]=0.\]

Assume that the Poisson cohomology $H^2(M,P_1)$ of the first Poisson structure is trivial, then any 
deformation $(\tilde{P}_1, \tilde{P}_2)$ of the bihamiltonian structure $(P_1, P_2)$ 
is equivalent to a deformation of the form
\begin{equation}\label{ky-4}
(\tilde{P}_1, \tilde{P}_2)=(P_1, P_2+\xi Q_{1}+\xi^2 Q_{2}+\dots),
\end{equation}
and $\tilde{P}_2$ must satisfy the equations
\begin{equation}
[P_1, \tilde{P}_2]=0,\quad  [\tilde{P}_2, \tilde{P}_2]=0.
\end{equation}
So the differential graded Lie algebra
$(\tilde{\Lambda}^*, [\ \,,\ ], d_2)$ controls the deformations of the bihamiltonian structure $(P_1, P_2)$.

A solution $Q_1\in \Lambda^2$ of the equations 
\begin{equation}
 [P_1, Q_1]=0, \quad [P_2, Q_1]=0
\end{equation}
gives a deformation 
\begin{equation}\label{ky-3}
(P_1, P_2)\to (P_1, P_2+\xi Q_1)
\end{equation}
of the bihamiltonian structure $(P_1, P_2)$ at the approximation up to 
the first order of $\xi$, and it is called an infinitesimal 
deformation of the bihamiltonian structure. We introduce the following cohomology groups
\begin{align}
&BH^0(M,P_1, P_2):=\left.\Ker d_2\right|_{\tilde{\Lambda}^0}=\Ker d_1|_{{\Lambda}^0}\cap \Ker d_2|_{{\Lambda}^0}\,,\label{zh-37}\\
&BH^k(M,P_1, P_2):=\frac{\left.\Ker d_2\right|_{\tilde{\Lambda}^k}}{\left.\Img d_2\right|_{\tilde{\Lambda}^{k-1}}}=\frac{\Ker d_1|_{{\Lambda}^k}\cap \Ker d_2|_{{\Lambda}^k} }{\Img d_2|_{\tilde{\Lambda}^{k-1}}}\,,\quad  k\ge 1,\label{zh-37-b}
\end{align}
which are called the bihamiltonian cohomologies of $(M, P_1, P_2)$. It is easy to see that
\begin{itemize}
\item[i)] The group $BH^2(M,P_1,P_2)$ characterizes the equivalence classes of infinitesimal deformations of the bihamiltonian structure $(P_1, P_2)$;
\item[ii)] If $BH^3(M,P_1,P_2)$ vanishes, then any infinitesimal deformation of the bihamiltonian structure $(P_1, P_2)$  can be extended to a deformation
\eqref{ky-4} of  $(P_1, P_2)$.
\end{itemize}

The above setting of Poisson structures and their cohomologies is generalized in \cite{DZ} to an infinite dimensional version. Poisson structures now live on the formal loop space $\mathcal{L}(M)$ of $M$, and they are given by Poisson bivectors on $\mathcal{L}(M)$, see also \cite{KKV, KKVV, Getz} and references in \cite{DZ} for similar formulations of infinite dimensional Poisson structures. 
Similar to the finite dimensional case, a Poisson bivector defines a Poisson bracket
on the space of local functionals of  $\mathcal{L}(M)$. This version of infinite dimensional Poisson structures is formulated for the purpose of 
studying systems of nonlinear evolutionary PDEs of the form
\begin{equation}\label{eq-hierarchy}
\frac{\p u^i}{\p t}=A^i_{k}(u)\,u^k_x+\e \left(B^i_{k}(u) u^k_{xx}+C^i_{k,l}(u) u^k_x u^l_x\right)+ \mathcal{O}(\e^2)
\end{equation}
which possess bihamiltonian structures with hydrodynamic limit. Here $\e$ is a dispersion parameter and $\mathcal{O}(\e^2)$ represents  the higher order in $\e$ terms $\sum_{k\ge 2}R^i_k \e^k$, and $R^i_k$ are polynomials in $u^{i,s}=\p_x^s u^i\ (s\ge 1)$ whose coefficients are smooth functions on $M$ (such polynomials are called differential polynomials). Here we set $\deg u^{i,s}=s$, and require that $\deg R^i_k=k+1$. In this way,
the power $k$ of $\e^k$ counts the degree of its coefficients.
We also assume summations for repeated upper and lower indices above and in what follows unless otherwise specified.
The system of evolutionary PDEs \eqref{eq-hierarchy} is said to possess a bihamiltonian structure with hydrodynamic limit if there 
exists a compatible pair of Poisson brackets of the form
\begin{align}\label{eq-bihamiltonian}
&\{u^i(x), u^j(y)\}_a=g^{ij}_a(u)\delta'(x-y)+\Gamma^{ij}_{k,a}(u)u^k_x\delta(x-y)\nn\\
&\quad + \sum_{k\ge 1}\e^k \sum_{l=0}^{k+1}P^{ij}_{kl,a}(u,u_x,\dots) \delta^{(k+1-l)}(x-y),\quad a=1,2
\end{align}
and Hamiltonians
\begin{equation}\label{ky-18}
H_a=\int_{S^1} h_a(u)dx +\sum_{k\ge 1}\e^k \int_{S^1}  f_{k,a}(u,u_x,\dots,u^{(k)}) dx, \quad  a=1,2 
\end{equation}
such that 
\begin{equation}\label{zh-20-a}
\frac{\p u^i}{\p t}=\{u^i(x), H_1\}_1=\{u^i(x), H_2\}_2.
\end{equation}
Here the matrices $(g^{ij}_a(u)),\ a=1,2$ are symmetric and nondegenerate, and $P^{ij}_{kl,a}, f_{k,a}$ are differential polynomials of degree $l$ and $k$ respectively. The integration $\int_{S^1}\cdot\,dx$ is taken on the class of functions
that are periodic in $x\in S^1$.

The bihamiltonian  recursion relation given by the second equality of \eqref{zh-20-a} enables one to find in general
infinitely many Hamiltonians which are in involution w.r.t. both of the Poisson brackets \cite{Magri-1, Getz, DZ, Magri-2}.
These Hamiltonians yield a hierarchy of bihamiltonian integrable systems which includes the originally given one \eqref{eq-hierarchy}. Many important nonlinear integrable hierarchies, including the ones associated to affine Lie algebras constructed by using the Drinfeld-Sokolov reduction procedure \cite{DS} and the ones that appear in Gromov-Witten theory and singularity theory (see \cite{Wit, Kon, Dij, Du1, DZ-cmp, DZ, Gi-2} and references therein),  are bihamiltonian systems of the above form. 

The leading terms of the r.h.s. of \eqref{eq-bihamiltonian} obtained by setting $\e=0$ 
give a compatible pair of Poisson brackets 
\begin{equation}\label{ky-5}
\{u^i(x), u^j(y)\}_a^{[0]}=g^{ij}_a(u)\delta'(x-y)+\Gamma^{ij}_{k,a}(u)u^k_x\delta(x-y),\quad  a=1,2
\end{equation}
of hydrodynamic type. Such a bihamiltonian structure is called semisimple if 
the characteristic polynomial $\det(g^{ij}_2-\lambda g^{ij}_1)$ in $\lambda$ has pairwise distinct and non-constant roots. The bihamiltonian structure \eqref{eq-bihamiltonian} now can be regarded as a deformation of \eqref{ky-5}.
In what follows, we always consider deformations of this form.
Unlike the finite dimensional case, here we use $\e$ as the deformation parameter to emphasize that
the power $k$ of $\e^k$ also count degrees of the deformation terms in the following sense:
we assume
$\deg \delta^{(s)}(x-y)=s$, then the coefficients of $\e^k$ in \eqref{eq-bihamiltonian} are homogeneous with
degree $k+1$.

A natural problem which is important in the theory of integrable systems and in its applications in mathematical physics is
to classify the bihamiltonian structures of the form \eqref{eq-bihamiltonian} and the associated integrable hierarchies. The
first step toward solving this problem is done in \cite{LZ1, DLZ1} with the help of the infinite dimensional version of the
bihamiltonian cohomologies formulated in \cite{DZ}.  Assuming the manifold $M$ is an open ball in $\R^n$,
it is proved in \cite{LZ1, DLZ1} that under the action of Miura type transformations
\begin{equation}\label{ky-30}
u^i\mapsto \tilde{u}^i=u^i+\sum_{k\ge 1} \e^k F^i_k(u,u_x,\dots,u^{(k)}),\quad i=1,\dots,n
\end{equation}
with $F^i_k$ being differential polynomials of degree $k$
(they are infinite dimensional analogue of 
gauge transformations, see Definition \ref{ky-7}, or its definitions given in \cite{DZ, jacobi}), any deformation of a semisimple bihamiltonian structure of hydrodynamic type \eqref{ky-5} can be transformed to a
bihamiltonian structure of the form \eqref{eq-bihamiltonian} which keeps the first Hamiltonian 
structure $\{u^i(x), u^j(y)\}^{[0]}_1$  undeformed and the deformation of the second Hamiltonian 
structure $\{u^i(x), u^j(y)\}^{[0]}_2$ contains only terms with even powers of $\e$ and, moreover,  two
deformations are equivalent under a Miura type transformation if and only if they are equivalent up to $\e^2$-order
approximation. So in what follows we will call 
a deformation of the bihamiltonian structure \eqref{ky-5} of the form
\begin{align}\label{ky-6}
&\{u^i(x), u^j(y)\}_1=\{u^i(x), u^j(y)\}_1^{[0]},\\
&\{u^i(x), u^j(y)\}_2=\{u^i(x), u^j(y)\}^{[0]}_2
+\e^2\sum_{l=0}^3 P^{ij}_{l}(u,u_x,u_{xx},u_{xxx}) \delta^{(3-l)}(x-y) \label{ky-6-b}
\end{align}
at the approximation up to $\e^2$-order an infinitesimal deformation. Here
$P^{ij}_{l}$ are differential polynomials of degree $l$.
It is also proved that any equivalent class (under the Miura type transformations) of infinitesimal deformations \eqref{ky-6}, \eqref{ky-6-b}
is characterized by a set of $n$ functions
$c_i(\lambda_i),\ i=1,\dots,n$ depending on the canonical coordinates $\lambda_i$ (see \cite{DZ} for their definition) of
the bihamiltonian structure. These functions are called the central invariants of the deformed bihamiltonian structure \eqref{eq-bihamiltonian}.
On the other hand, any set of smooth functions $c_1(\lambda_1),\dots, c_n(\lambda_n)$ determines an infinitesimal deformation \eqref{ky-6}, \eqref{ky-6-b}
of a given semisimple bihamiltonian structure of hydrodynamic type \eqref{ky-5}.

In the present and the subsequent paper \cite{DLZ2} we are to consider the problem of whether one can extend any infinitesimal
deformation \eqref{ky-6}, \eqref{ky-6-b} of a semisimple bihamiltonian structure of hydrodynamic type \eqref{ky-5} to a deformation of the form \eqref{eq-bihamiltonian}. To this end, we
reformulate in this paper the notion of infinite dimensional Poisson structures in terms of the infinite jet space 
$J^\infty(\hat{M})$ of a super manifold $\hat{M}$ (see Section 2 for definition). This reformulation provides us with a more
convenient way to study properties of the bihamiltonian cohomologies for a semisimple bihamiltonian structure of hydrodynamic
type. In particular, the long exact sequence \eqref{bhc-long} which we prove in Corollary \ref{cor-44} provides an important
tool to compute the bihamiltonian cohomologies. By using this result, we compute in Section \ref{sec-5} the third bihamiltonian
cohomology
\[BH^3(\hat\F)=\bigoplus_{d\ge 0} BH^3_d(\hat\F)\]
of the following bihamiltonian structure
\begin{align}
&\{u(x),u(y)\}_1=\delta'(x-y), \label{zh-19}\\
&\{u(x),u(y)\}_2=u(x)\delta'(x-y)+\frac12 u_x(x) \delta(x-y) \label{zh-20}
\end{align}
for the dispersionless KdV hierarchy
\begin{equation}\label{zh-31}
\frac{\p u}{\p t_p}=\frac1{p!}u^p u_x,\quad  p\ge 0.
\end{equation}
The equation \eqref{zh-31} of this hierarchy can be represented as a bihamiltonian system of  the form \eqref{zh-20-a} with
\[H_{1,p}=\frac1{(p+2)!} \int_{S^1} u(x)^{p+2} dx,\quad  H_{2,p}= \frac{2}{(2p+1)(p+1)!}  \int_{S^1} u(x)^{p+1} dx.\]
We prove the following theorem in the present paper:
\begin{thm}\label{thm-zh1}
For the bihamiltonian structure \eqref{zh-19}, \eqref{zh-20} the bihamiltonian cohomologies $BH^3_d(\hat\F)$ are trivial for $d\ge 4$.
\end{thm}

For the bihamiltonian structure \eqref{zh-19}, \eqref{zh-20} the canonical coordinate is given by $\lambda=u$. As a corollary of the above theorem, we have
\begin{thm}\label{thm-zh2}
For any given smooth function $c(u)$ the bihamiltonian structure \eqref{zh-19}, \eqref{zh-20} has a unique equivalent class of deformations with a representative of the form
 \begin{align}
&\{u(x),u(y)\}_1=\delta'(x-y), \label{zh-23}\\
&\{u(x),u(y)\}_2=u(x)\delta'(x-y)+\frac12 u(x) \delta(x-y)\nn\\
&\quad +\e^2 \left(3 c(u) \delta'''(x-y)+\frac92 c'(u) u_x \delta''(x-y) +\frac32\,c''(u) u_x^2 \delta'(x-y)\right.\nn\\
&\quad\quad \left.+\frac32 c'(u) u_{xx}\delta'(x-y)\right)+\sum_{g\ge 2}\e^{2 g}\sum_{k=0}^{2g+1} A_{g,k}(u,u_x,\dots,u^{(k)})
\delta^{(2g+1-k)}(x-y).\label{zh-24}
\end{align}
Here $A_{g;k}$ are differential polynomials of degree $k$.
\end{thm} 
The validity of this theorem was conjectured by Lorenzoni in \cite{Loren}. He also gave
the approximation of the deformation up to $\e^4$. In \cite{AL}, Arsie and Lorenzoni further extend the
approximation up to $\e^8$.

The bihamiltonian structures in the above theorem for the particular choices of $c(u)=\frac1{24}$ and $c(u)=\frac{u}{24}$
are well known in the theory of integrable systems. They provide bihamiltonian structures for the KdV hierarchy and the
Camassa-Holm hierarchy respectively \cite{LZ1}. For other choices of the central invariant $c(u)$ the above bihamiltonian
structures and the associated integrable hierarchies are new. In the present paper we will consider in some detail the
bihamiltonian structure and the associated integrable hierarchy when $c(u)$ is inversely proportional to $u$.

The paper is organized as follows. In Section \ref{sec-2} we define the differential polynomial algebra and the Schouten bracket on the infinite jet space of a super manifold. In Section \ref{sec-3} we formulate the infinite dimensional Poisson structure and their cohomologies in terms of the notions established in Section \ref{sec-2}. In Section \ref{sec-4} we define the notion of bihamiltonian cohomologies, and present an important method
to compute the bihamiltonian cohomologies. In Section \ref{sec-5} we compute the third bihamiltonian cohomology of the    
bihamiltonian structure \eqref{zh-19}, \eqref{zh-20}, and prove Theorem \ref{thm-zh1} and Theorem \ref{thm-zh2}. In Section \ref{sec-6} we consider  examples of the bihamiltonian structures corresponding to some particular choices of the central invariant $c(u)$. The final section is for conclusion.

\section{The differential polynomial algebra and the Schouten bracket} \label{sec-2}
Let $M$ be a smooth manifold of dimension $n$, $\hM$ be the super manifold $\Pi(T^*M)$, i.e. the cotangent bundle of $M$
with fiber's parity reversed.
It is well known that the algebra of multi-vectors, i.e. sections of the exterior algebra  bundle $\Lambda(TM)$
of the tangent bundle $TM$, can be regarded as the algebra of smooth functions on $\hM$, that is
\[C^\infty(\hat{M})=\Gamma(\Lambda(TM)).\]
With the help of this identification, the Schouten-Nijenhuis bracket between multi-vectors can be computed via the Poisson bracket of
the canonical (super) symplectic structure on $\hM$. Suppose $P\in \Gamma(\Lambda^p(TM))$, $Q\in \Gamma(\Lambda^q(TM))$
are two multi-vectors with degree $p$ and $q$ respectively, their Schouten-Nijenhuis bracket $[P, Q]$, which is a multi-vector of degree $p+q-1$,
has the following expression on a local trivialization $\hU=U\times \R^{0|n}$ of $\hM$
\begin{equation}
[P,Q]=\frac{\p P}{\p \theta_i}\frac{\p Q}{\p u^i}+(-1)^p\frac{\p P}{\p u^i}\frac{\p Q}{\p \theta_i}, \label{bra-fd}
\end{equation}
where $u^1,\dots, u^n$ are coordinates on $U$, and $\theta_1,\dots, \theta_n$ are the corresponding dual coordinates on the fiber $\R^{0|n}$.

In this section, we introduce the infinite dimensional analog of \eqref{bra-fd} on the infinite jet space $J^{\infty}(\hM)$.

Let $J^k(\hM)\ (k\ge 1)$ be the $k$-th jet space of $\hM$. The manifold $J^k(\hM)$ is, by definition, a fiber bundle with fibers being the spaces of $k$-th
Taylor polynomials of germs of curves on $\hM$. Let $\hU$ be a local chart of $\hM$ with local coordinates $(u^i,  \theta_i)\ (i=1, \dots, n)$,
and $\hat{U}\times \R^{nk|nk}$ be a local trivialization of $J^k(\hM)$ over $\hU$. We can take the coordinates on the fiber $\R^{nk|nk}$
as the values of higher derivatives of curves
\[u^{i,s}=\left.\frac{d^s}{dx^s}u^i(x)\right|_{x=0},\quad \theta_i^s=\left.\frac{d^s}{dx^s}\theta_i(x)\right|_{x=0},
\quad s=1,\dots,k.\]
It is easy to see that the transition functions of $J^k(\hM)$, induced by the change of coordinates on $\hM$
from $(u^i, \theta_i)$ to $(\tilde{u}^i, \tilde{\theta}_i)$ with
\[\tilde{u}^i=\tilde{u}^i(u^1,\dots,u^n),\quad \tilde{\theta}_i=\frac{\p u^j}{\p \tilde{u}^i}\theta_j\,, \]
are given by the chain rule of higher derivatives as follows:
\begin{equation}\label{tran-2}
\tilde{u}^{i,s}=\sum_{t=0}^{s-1}\frac{\p \tilde{u}^{i,s-1}}{\p u^{j,t}}u^{j,t+1},\quad
\tilde{\theta}_i^s=\sum_{t=0}^s\binom{s}{t}\frac{\p u^{j,t}}{\p \tilde{u}^i}\theta_j^{s-t}. 
\end{equation}
For example, when $s=1, 2$ we have
 \begin{align*}
\tilde{u}^{i,1}=&\frac{\p \tilde{u}^i}{\p u^j}u^{j,1},\\
\tilde{\theta}_i^1=&\frac{\p \tilde{\theta}_i}{\p \theta_j}\theta_j^1+\frac{\p \tilde{\theta}_i}{\p u^j}u^{j,1}
=\frac{\p u^j}{\p \tilde{u}^i}\theta_{j}^1+\frac{\p^2 u^j}{\p \tilde{u}^i\p \tilde{u}^l}\tilde{u}^{l,1}\theta_j,\\
\tilde{u}^{i,2}=&\frac{\p \tilde{u}^i}{\p u^j}u^{j,2}+\frac{\p^2 \tilde{u}^i}{\p u^l \p u^m}u^{l,1}u^{m,1}\\
\tilde{\theta}_i^2=&\frac{\p u^j}{\p \tilde{u}^i}\theta_j^2+2\frac{\p^2 u^j}{\p \tilde{u}^i\p \tilde{u}^l}\tilde{u}^{l,1}\theta_j^1\\
&+\left(\frac{\p u^j}{\p^2 \tilde{u}^i\p \tilde{u}^l}\tilde{u}^{l,2}+\frac{\p^3 u^j}{\p \tilde{u}^i\p \tilde{u}^l\p \tilde{u}^m}
\tilde{u}^{l,1}\tilde{u}^{m,1}\right)\theta_j.
\end{align*}
Note that these transition functions are always polynomials of the jet variables $u^{i,s}, \theta_i^s$
$(s\ge1)$
with coefficients being smooth functions on $M$.

For $k \ge l$, there exists a natural forgetful map
\[\pi_{k,l}:J^k(\hM)\to J^l(\hM),\]
which just forgets the derivatives with orders greater than $l$. It is easy to see that the system
\[\left(\{J^k(\hM)\}_{k\ge1},\ \{\pi_{k,l}\}_{k\ge l\ge1}\right)\]
forms a projective system, so we can take its projective limit
\[J^\infty(\hM)=\varprojlim_k J^k(\hM),\]
which is called the infinite jet space of $\hM$.

The algebras of smooth functions on $\{J^k(\hM)\}_{k\ge1}$ and the pullback maps
$\{\pi^*_{k,l}\}_{k\ge l\ge1}$ among them form an inductive system, so we can define its inductive limit
\[C^\infty(J^\infty(\hM))=\varinjlim_k C^\infty(J^k(\hM)),\]
which is called the algebra of smooth functions on $J^\infty(\hM)$.

A function $f\in C^\infty(J^\infty(\hM))$ is called a differential polynomial if it depends on the jet variables polynomially in
certain local coordinate system. Due to the form of transition functions  \eqref{tran-2}, this definition is independent of
the choice of local coordinate system. All differential polynomials form a subalgebra of $C^\infty(J^\infty(\hM))$, we denote this
subalgebra by $\AA$. Locally, we can regard $\AA$ as
\[C^\infty(\hU)[u^{i,s}, \theta_i^s\mid i=1, \dots, n,\ s=1, 2, \dots],\]
where $\hU$ is a chart on $\hM$ with local coordinates $(u^i, \theta_i)\ (i=1, \dots,n)$.

The space $\AA$ is not big enough to define, say, Miura type transformations and their inverses, so we need to enlarge it.
First we introduce a gradation on $\AA$ by defining
\[\deg u^{i,s}=\deg \theta_i^s=s \mbox{ for } s>0,\mbox{ and } \deg f=0 \mbox{ for } f\in C^\infty(\hM),\]
then for any element $f\in\AA$, one can decompose it into a direct sum of homogeneous components
\[f=f_0+f_1+\dots, \quad \mbox{where }  \deg f_k=k.\]
This gradation induces the valuation
\begin{equation}\label{ky-9}
\nu(f)=\left\{\begin{array}{ll}\min\{k\mid f_k\ne0\}, & \quad f\ne0 \\ \infty, & \quad f=0 \end{array}\right. 
\end{equation}
We denote the distance induced by this valuation by
\[d(f,g)=e^{-\nu(f-g)}\]
and complete $\AA$ by using this distance.
By abuse notations, we still denote the completion of $\AA$ by $\AA$, and still call it the differential polynomial algebra of $\hM$.
Locally, the completed $\AA$ looks like
\[C^\infty(\hU)[[u^{i,s}, \theta_i^s\mid i=1, \dots, n,\ s=1, 2, \dots]].\]
Note that the topology on $C^\infty(\hU)$ is, by definition, the discrete one, so only finite sum of smooth functions are allowed.
Elements of $\AA$ are still called differential
polynomials, though they may be formal series of homogeneous differential polynomials with strict increasing degrees.
We call the above defined  gradation the \textit{standard gradation}, and denote the degree $d$ component of $\AA$ w.r.t this gradation by $\AA_d$.

Sometimes we use a formal parameter $\e$ to indicate the degree of homogeneous terms, i.e. we represent $f\in\AA$ in the form
\begin{equation}\label{ky-19}
f=f_0+\e\,f_1+\e^2\,f_2+\dots,
\end{equation}
where $f_k\in\AA_{k+d_0}$, and $d_0=\nu(f)$.
In  this way, the topology on $\AA$ is just the $\e$-adic topology.

The super variables in $\AA$ induce another gradation
\[\deg \theta_i^s=1,\quad \deg u^{i,s}=\deg f=0.\]
We call it \textit{super gradation}, and denote the degree $p$ component of $\AA$ w.r.t this gradation by $\AA^p$. We also use the notation
$\AA^p_d=\AA^p\cap \AA_d$. The degree $0$ component $\AA^0$ is often denoted by $\A$, and is called the differential polynomial algebra
on $M$. It's easy to see that $\AA^0_0=C^\infty(M)$.

By using the form of the transition functions given in \eqref{tran-2}, one can show that the vector field
\[\p=\sum_{s\ge0}\left(u^{i,s+1}\frac{\p}{\p u^{i,s}}+\theta_i^{s+1}\frac{\p}{\p \theta_i^s}\right)\]
is globally defined on $J^\infty(\hM)$, where $u^{i,0}=u^i$ and $\theta_i^0=\theta_i$.
It defines a derivation on $\AA$, whose restriction on $\A$ also yields a derivation on $\A$,
so we have the following definition.

\begin{dfn}
We denote the quotient space $\AA/\p\,\AA$ by $\FF$ and call its elements local functionals on $\hM$. 
We also denote the subspace $\FF^0=\A/\p\,\A$ of $\FF$ by $\F$ and call its elements local functional on $M$. 
\end{dfn}
Note that $\p$ preserves both the standard gradation and the super gradations
on $\AA$, so the quotient space $\FF$ also possesses the two gradations. We
also use the notations $\FF^p$, $\FF_d$, and $\FF^p_d$ to represent the
corresponding homogeneous components of $\FF$ with respect to these gradations.
Elements in $\FF^p$ are often called $p$-vectors. In particular, for the $p=1, 2, 3$ cases,
we will use the names vector, bivector, and trivector respectively.

In what follows we will use $\int$ to denote the projection $\A\to\F$ and $\AA\to\FF$.

\begin{dfn}
We denote by $\E=\Der(\A)^\p$ the centralizer of $\p$ in the Lie algebra of continuous derivations of $\A$ and call its elements evolutionary vector fields on $M$.
\end{dfn}

From the above definition it follows that any evolutionary vector field takes the form
\begin{equation}\label{zh-34}
D_X=\sum_{s\ge0}\p^s(X^i)\frac{\p}{\p u^{i,s}},\quad \mbox{where } X^i\in\A.
\end{equation}
If we define
\begin{equation}
X=\int X^i\theta_i\,, \label{eq-intx}
\end{equation}
then one can verify that $X$ is globally defined on $\hM$, so it is in fact an element of $\FF^1$.
On the other hand, each element of $\FF^1$ has a representative of the form \eqref{eq-intx}, so it corresponds to an
element $D_X\in\E$.  Thus  we will identify $\E$ with $\FF^1$ in this way from now on.

We can also associate to any element \eqref{zh-34} of $\E$  the following system of evolutionary PDEs 
\begin{equation}\label{zh-32}
u^i_t=X^i.
\end{equation}
If $u=(u^1,\dots,u^n)$ is a solution of this system, and $f$ is a differential polynomial in $\A$, then we have
\[f_t=D_X(f).\]
So a local functional $F\in \F$ is a conserved quantity of \eqref{zh-32} if and only if  $D_X(F)=0$.
Here we used the fact that $\F$ is a module of the Lie algebra $\E$.

Now let us define an infinite dimensional analog of the bracket \eqref{bra-fd} on the space $\FF$.
To this end we first introduce the variational derivatives of $f\in\AA$ w.r.t $\theta_i$ and $u^i$ as follows:
\[\delta^i f=\sum_{s\ge0}(-\p)^s\frac{\p f}{\p \theta_i^s},\quad  \delta_i f=\sum_{s\ge0}(-\p)^s\frac{\p f}{\p u^{i,s}},
\quad  i=1,\dots,n.\]
It is well known that
\[\delta^i\,\p=0,\quad \delta_i\,\p=0,\]
so $\delta^i$ and $\delta_i$
induce maps from $\FF$ to $\AA$. We denote the induced maps by the same notations $\delta^i$ and $\delta_i$.
Then we can define a bracket
\[[\ \,,\ ]:\FF^p\times \FF^q \to \FF^{p+q-1}, \quad (P,Q)\mapsto [P, Q]\]
with
\begin{equation}
[P,Q]=\int\left(\delta^i P\,\delta_i Q+(-1)^p \delta_i P\,\delta^i Q\right). \label{bra-id}
\end{equation}

\begin{thm}[\cite{jacobi}]\label{thm-bra}
For any $P\in\FF^p$, $Q\in\FF^q$, $R\in\FF^r$, we have
\begin{itemize}
\item[i)] $[P, Q]=(-1)^{p\,q}[Q,P]$;
\item[ii)] $(-1)^{p\,r}[[P, Q],R]+(-1)^{q\,p}[[Q, R],P]+(-1)^{r\,q}[[R, P],Q]=0$.
\end{itemize}
\end{thm}
The bracket $[\ \,,\ ]$ is called the Schouten-Nijenhuis bracket on $J^\infty(\hM)$.
From now on, we call it the Schouten bracket for short.

In order to compute the bihamiltonian cohomologies, we also need to define another useful map.
Let us denote by $\EE^p_d$ the space of super derivations of degree $(p,d)$ over $\AA$.
Here a super derivation of degree $(p,d)$ over $\AA$ is a continuous linear map $\Delta:\AA \to \AA$
satisfying $\Delta(\AA_{d'}^{p'})\subset\AA_{d'+d}^{p'+p}$, and 
\[\Delta(f\cdot g)=\Delta(f)\cdot g+(-1)^{p\,k}f\cdot \Delta(g),  \quad \mbox{where } f\in\AA^k.\]
We also denote
\[\EE_d=\prod_{p\in\Z}\EE^p_d,\quad  \EE^p=\prod_{d\in\Z}\EE^p_d,\quad \EE=\prod_{d\in\Z}\prod_{p\in\Z}\EE^p_d.\]
We define the following map
\[D: \FF^p \to \EE^{p-1},\quad P\mapsto D_P\]
with
\begin{equation}\label{zh-35}
D_P=\sum_{s\ge0}\left(\p^s(\delta^i P)\frac{\p}{\p u^{i,s}}+(-1)^p\p^s(\delta_i P)\frac{\p}{\p \theta_i^s}\right).
\end{equation}
Note that for $P=X\in\FF^1$ the operator $D_P$, when restricted to $\A$, coincides with the one that is defined in \eqref{zh-34}.

It is well known that the space $\EE$ has a Lie superalgebra structure  w.r.t. the super degree, i.e. for $D_1\in\EE^{p_1}_{d_1}$,
$D_2\in\EE^{p_2}_{d_2}$ the bracket
\[[D_1, D_2]=D_1\,D_2-(-1)^{p_1\,p_2}D_2\,D_1\in\EE^{p_1+p_2}_{d_1+d_2}\]
is a super Lie bracket over $\EE$.

Now let us give some examples to illustrate the notions that we introduced above. We assume $\dim M=1$, and denote $u^{i, s}, \theta_i^s$ by $u^s$ and $\theta^s$ respectively. Take 
\[P= \int u u^1\theta\in\FF^1_1,\quad Q=\int u^2 \theta\theta^1\in\FF^2_3.\]
Then from the definition \eqref{bra-id} of the Schouten bracket we have
\begin{align*}
[P,Q]&=\int \left(u u^1 \p^2(\theta\theta^1)-(u^1\theta-\p (u\theta))(u^2 \theta^1+\p(u^2 \theta))\right)\\
&=\int\left(u u^1 (\theta^1\theta^2+\theta\theta^3)+u \theta^1 (2 u^2 \theta^1+u^3\theta)\right)\\
&=\int\left(u u^1\theta^1\theta^2+u u^1\theta\theta^3-u u^3\theta\theta^1\right)\in\FF^2_4,
\end{align*}
and the operators $D_P, D_Q$ are given by
\begin{align*}
&D_P=\sum_{s\ge 0} \left(\p^s (u u^1)\frac{\p}{\p u^s}+\p^s (u \theta^1) \frac{\p}{\p\theta^s}\right)\in\EE^0_1,\\
&D_Q=\sum_{s\ge 0} \left(\p^s(2 u^2\theta^1+u^3\theta)\frac{\p}{\p u^s}+\p^{s+2}(\theta\theta^1)\frac{\p}{\p \theta^s}\right)\in\EE^1_3.
\end{align*}
These operators satisfy the identity $D_{[P,Q]}=[D_P, D_Q]$. In general we have
the following theorem.
\begin{thm}[\cite{jacobi}]\label{thm-bra2}
For any $P\in\FF^p$, $Q\in\FF^q$ the following identities hold true:
\begin{itemize}
\item[i)] $[P,Q]=\int D_P(Q)$;
\item[ii)] $D_{[P,Q]}=(-1)^{p-1}[D_P, D_Q]$.
\end{itemize}
\end{thm}

\begin{rmk}
If we regard $P\in\FF^p$ as a Hamiltonian, then $D_P$ is the corresponding Hamiltonian vector field, and the identity \textit{ii)} of the above
theorem is just the analog of the homomorphism from the Poisson algebra of a symplectic manifold to the Lie algebra of its Hamiltonian vector
fields.
\end{rmk}

\section{Hamiltonian structures and their cohomologies} \label{sec-3}

Now we are ready to formulate the infinite dimensional version of Poisson structures which we also call Hamiltonian structures.
We assume from now on that the manifold $M$ is an open ball in $\R^n$.

\begin{dfn}
A bivector $P\in\FF^2$ is called a Poisson bivector or a Hamiltonian structure if $[P, P]=0$. We say that an evolutionary vector field $X\in \E$ has a Hamiltonian structure if there exists a Poisson bivector $P\in \FF^2$ and a local functional
$F\in\F$ such that $X=-[P,F]$.
\end{dfn}

For a given Poisson bivector $P$, we define the associated Poisson bracket as follows:
\[\{\ ,\ \}_P:\F\times\F\to\F,\quad (F, G)\mapsto\{F, G\}_P,\]
where $\{F, G\}_P=-[F, [P, G]]$. Note that for any $P\in\FF^2$, we can uniquely find $P^{ij}_s\in\A$ such that
\[P=\frac12\int P^{ij}_s\theta_i\theta_j^s,\]
and
\begin{equation}\label{zh-80}
P^{ij}_s \p^s=(-1)^{s+1}\p^s P^{ji}_s.
\end{equation}
Then the Poisson bracket associated to $P$ has the form:
\begin{equation}\label{poisson}
\{F, G\}_P=\int\frac{\delta F}{\delta u^i}\left(P^{ij}_s \p^s\right)\left(\frac{\delta G}{\delta u^j}\right).
\end{equation}
The operator $\hat{P}=P^{ij}_s\p^s$ is called the Hamiltonian operator of $P$. In literatures
on integrable systems, the Hamiltonian structures of partial differential
equations are often given by such operators.

\begin{emp}
Suppose $P\in \FF^2_0$, then $P$ must take the following form
\[P=\frac12\int P^{ij}(u)\theta_i\,\theta_j, \quad \mbox{with }  P^{ij}(u)=-P^{ji}(u),\]
so it actually corresponds to the following bivector on $M$:
\[\bar{P}=\frac12P^{ij}(u)\frac{\p}{\p u^i}\wedge\frac{\p}{\p u^j}.\]
It's easy to see that the bracket $[P,P]$ (see \eqref{bra-id}) coincides with $[\bar{P},\bar{P}]$ (see \eqref{bra-fd}),
so $P$ is a Poisson bivector if and only if $\bar{P}$ is a Poisson bivector on $M$.
\end{emp}

\begin{emp}[\cite{DN83}]
Suppose $P\in \FF^2_1$, then $P$ can be represented as 
\[P=\frac12\int\left(g^{ij}(u)\theta_i\theta_j^1+\Gamma^{ij}_k(u)u^{k,1}\theta_i\theta_j\right) \]
which satisfies the normalization condition \eqref{zh-80}.
We also assume that the matrix $(g^{ij})$ is nondegenerate, then $P$ is a Hamiltonian structure if and only if
\begin{itemize}
\item [i)] $g=(g_{ij})=(g^{ij})^{-1}$  is a flat metric on $M$;
\item [ii)] $\Gamma_{lk}^j=-g_{li}\Gamma^{ij}_k$ are the Christoffel symbols of the Levi-Civita connection of $g$.
\end{itemize}
Hamiltonian structures of this form are called of hydrodynamic type.
\end{emp}

\begin{rmk} \label{ky-8}
In this remark, we explain the relationship between the notation of Hamiltonian structures introduced via the loop space
approach (see Section 1, Section 6, and \cite{DZ}) and the one via the jet space approach used here.

Let $S^1$ be the one-dimensional circle, $\CL(M)=C^\infty (S^1, M)$ be the free loop space of $M$.
According to the definition of $J^\infty(M)$, every loop $\phi \in \CL(M)$ can be lifted to a loop
$\phi^\infty\in \CL(J^\infty (M))=C^\infty(S^1, J^\infty(M))$.
For a differential polynomial 
$f\in \A$  on $J^\infty(M)$, one can define a
functional on $\CL(M)$ with values in $\R[[\e]]$ as follow:
\[\hat{f}[\phi]=\int_{S^1} \left(\phi^\infty\right)^* (f) dx=\int_{S^1} f(\phi^\infty(x))dx.\]
One can show that $\hat{f}$ vanishes if and only if $f$ is in the image of $\p$, so the linear space of functionals
of the above form is isomorphic to the space of local functionals $\F=\A/\p\A$.

We can also consider the functionals in the sense of distributions:
\[\phi \mapsto \int_{S^1} f(\phi^\infty(x))\delta^{(s_1)}(x-y_1)\cdots\delta^{(s_k)}(x-y_k)dx,\]
where $y_1, \dots, y_k \in S^1$. We denote by $\F'$ the linear space of functionals of this form, then the Poisson
bracket \eqref{poisson} can be extended to $\F'$ in the following way:
\[\{F, G\}_P[\phi]=\int_{S^1} \frac{\delta F}{\delta u^i(x)}
\left(P^{ij}_s(\phi^\infty(x))\p_x^s\right)\frac{\delta G}{\delta u^j(x)}dx.\]
Note that the result of the bracket is also a functional in $\F'$.

In particular, if we take
\[u^i(y)=\int_{S^1} u^i(x)\delta(x-y)dx,\quad u^j(z)=\int_{S^1} u^j(x)\delta(x-z)dx,\]
then
\begin{align*}
\{u^i(y), u^j(z)\}_P=&\int_{S^1} \delta(x-y)\left(P^{ij}_s(\phi^\infty(x))\p_x^s\right)\delta(x-z)dx\\
=&P^{ij}_s(u(y), u'(y), \dots)\delta^{(s)}(y-z).
\end{align*}
This is just the notation used in Section 1, Section 6 and \cite{DZ}.
For example, the Poisson bivectors $P_a,\, a=1,2$ corresponding to the bihamiltonian structure \eqref{eq-bihamiltonian} are given by
\begin{align*}
&P_a=\frac12 \int \left(g^{ij}_a(u)\theta_i \theta_j^1 +\Gamma^{ij}_{k,a}(u)u^k_x\theta_i \theta_j 
+ \sum_{k\ge 1}\e^k \sum_{l=0}^{k+1}P^{ij}_{kl,a}(u,u_x,\dots) \theta_i \theta_j^{k+1-l}\right),
\end{align*}
and the elements of $\FF^0$ that correspond to the Hamiltonians given in \eqref{ky-18} has the expressions
\[H_a=\int \left( h_a(u)+\sum_{k\ge 1}\e^k f_{k,a}(u,u_x,\dots,u^{(k)}) \right), \quad  a=1,2. \]
For the insertion of the parameter $\e$ in the expressions of $P_a$ and $H_a$
see \eqref{ky-19} for explanation.
\end{rmk}

By using Theorem \ref{thm-bra}, \ref{thm-bra2} and the definition of Hamiltonian structures, one can easily prove the following lemma.
\begin{lem}
Let $P\in\FF^2$ be a Hamiltonian structure, $d_P$ be its adjoint action, i.e.
\[d_P(Q)=[P,Q],\quad \forall\, Q\in\FF.\]
Then we have
\[d_P^2=0,\quad D_P^2=0.\]
\end{lem}
This lemma shows that both $(\FF, d_P)$ and $(\AA, D_P)$ are complexes. In fact, they are also DGLAs.
The graded Lie bracket on $\FF$ is just the Schouten bracket. The bracket on $\AA$, which will not be
used in the present paper, can be found in \cite{Getz}. The DGLA $(\FF, [\ \,,\ ], d_P)$ controls the 
deformation theory of the Hamiltonian structure $P$.

\begin{dfn}\label{ky-7}
\mbox{}
\begin{itemize}
\item[i)] Let $P\in\FF^2$ be a Hamiltonian structure,
a bivector $\tilde{P}=P+Q\in\FF^2$ is called a 
deformation of $P$ if $\nu(Q)>\nu(P)$ and $[\tilde{P}, \tilde{P}]=0$ (see \eqref{ky-9} for the definition of $\nu$).
\item[ii)] For two deformations $\tilde{P}_1, \tilde{P}_2$ of the Hamiltonian structure  $P$, if there exists $X\in\FF^1$ such that $\nu(X)>0$ and
\[\tilde{P}_2=\exp (\mathrm{ad}_X) \tilde{P}_1,\]
then we say that $\tilde{P}_2$ is equivalent to $\tilde{P}_1$.
\item[iii)] Let $P\in\FF^2$ be a Hamiltonian structure, an infinitesimal deformation of $P$ is a bivector $Q\in\FF^2$
such that $\nu(Q)>\nu(P)$ and $[P, Q]=0$.
\item[iv)] For two infinitesimal deformation $Q_1, Q_2$ of the Hamiltonian structure $P$, if there exists $X\in\FF^1$ such that $Q_2=Q_1+[P, X]$, then we say that $Q_2$ is equivalent to $Q_1$.
\end{itemize}
\end{dfn}
Note that if we rewrite $P, Q$ in the loop space notations (see Remark \ref{ky-8}), and insert powers of $\e$
to indicate the degree, then the above definition of deformation coincides with the one given in Section 1.
The gauge transformation appearing in the above equation is also called a Miura type transformation. They correspond to
coordinates transformations \eqref{ky-30} on the jet space $J^\infty(M)$, see
\cite{jacobi} for more details.

We are interested in the deformation problem of Hamiltonian structures of hydrodynamic type. In this case,
the Hamiltonian structure $P$ is homogeneous of degree 1, i.e. $P\in\FF^2_1$.
If the Hamiltonian structure $P$ is homogeneous of degree $d_0$,
then the cohomology group $H^k(\FF, d_P)$ has a gradation induced from the standard gradation of $\FF$:
\[H^k(\FF, d_P)=\bigoplus_{d\ge0}H^k_d(\FF, d_P).\]
Here
\[  
H^k_d(\FF, d_P)=\frac{\FF^k_d\cap\Ker d_P}{\FF^k_d\cap\Img  d_P},\quad k, d\ge 0.\]
It's easy to see that
\begin{itemize}
\item[i)] The group $H^2_{\ge d_0+1}(\FF, d_P)$ classifies the equivalence classes of infinitesimal deformations of
$P$;
\item[ii)] If $H^3_{\ge 2\delta}(\FF, d_P)$ vanishes, then any infinitesimal deformation
of $P$ can be extended to a deformation of it. Here $\delta$ is the lowest degree of classes in $H^2_{\ge d_0+1}(\FF)$.
\end{itemize}

The computation of $H^*(\FF,d_P)$ is not easy, so we introduce the complex $(\AA, D_P)$ to help us. Let $\R$ be the
subalgebra of $\AA$ consists of constant functions on $J^\infty(\hM)$, then one can show that $\R$ is the kernel
of $\p:\AA\to\AA$. Note that $\R$ is contained in $D_P$'s kernel, so we can define another complex $(\AA/\R, D_P)$, and have the following lemma.

\begin{lem}
The following sequence of complexes is exact
\[0\xrightarrow{\ \ \ }(\AA/\R, D_P)\xrightarrow{\ \p\ }(\AA, D_P)\xrightarrow{\ \int\ }(\FF, d_P)\xrightarrow{\ \ \ }0,\]
hence we have a long exact sequence of cohomologies,
\begin{equation}\label{long-ex-sq}
\cdots\xrightarrow{\ \ \ }H^p(\AA/\R)\xrightarrow{\ \ \ }H^p(\AA)\xrightarrow{\ \ \ }
H^p(\FF)\xrightarrow{\ \ \ }H^{p+1}(\AA/\R)\xrightarrow{\ \ \ }\cdots.
\end{equation}
\end{lem}

Now we assume that  $P$ is a Hamiltonian structure of hydrodynamic type, so $P\in\FF^2_1$ and $D_P\in \EE^1_1$.
Then the homologies of the above three complexes become direct sums of their homogeneous components, and the long
exact sequence \eqref{long-ex-sq} can be written as long exact sequences of homogeneous components:
\begin{equation}\label{long-ex-sq-homo}
\cdots\xrightarrow{\ \ \ }H^p_{d-1}(\AA/\R)\xrightarrow{\ \ \ }H^p_d(\AA)\xrightarrow{\ \ \ }H^p_d(\FF)
\xrightarrow{\ \ \ }H^{p+1}_d(\AA/\R)\xrightarrow{\ \ \ }\cdots.
\end{equation}
Note that $\R\subset \AA^0_0$, so we have
\[H^p_d(\AA/\R, D_P)=\left\{\begin{array}{ll}
H^p_d(\AA, D_P), &\quad  (p,d)\ne (0,0);\\
H^0_0(\AA, D_P)/\R,&\quad  (p,d)=(0,0).
\end{array}\right.\]

\begin{lem}[\cite{jacobi}] \label{lem-36}
Let $P\in\FF^2_1$ be a Hamiltonian structure of hydrodynamic type, then for $p\ge 0$ we have
\begin{equation}
H^p_d(\AA, D_P)=\left\{
\begin{array}{ll}0, &\quad d>0;\\ \bigwedge^p(\R^n), &\quad d=0.\end{array}\right.\end{equation}
\end{lem}

From the long exact sequence \eqref{long-ex-sq-homo} and the above lemma we 
obtain the following theorem \cite{jacobi}:
\begin{thm}\label{zh-thm-36}
Let $P\in\FF^2_1$ be a Hamiltonian structure of hydrodynamic type, then for $p\ge 0$ we have
\begin{equation}
H^p_{>0}(\FF, d_P)=0.
\end{equation}
\end{thm}
This theorem is first proved in \cite{Getz}, and it is also proved in \cite{Magri-2, DZ} for the $p=1,\ 2$ cases.
Recently, De Sole and Kac \cite{DSK1, DSK2} prove some theorems on cohomologies of Poisson vertex algebra,
which also imply the above theorem. In particular, the fact that $H^2_{>1}(\FF, d_P)=0$ implies that any
deformation of $P$ can be eliminated by a Miura type transformation (see \cite{LZ1, jacobi}).

\section{Bihamiltonian structures and their cohomologies}\label{sec-4}

In this section we give the definition of bihamiltonian cohomologies and show its main properties.

\begin{dfn}
Let $P_1, P_2\in \FF^2$ be two Poisson bivectors, if $[P_1, P_2]=0$ then we call the pair $(P_1, P_2)$ a bihamiltonian structure.
\end{dfn}

Let $(P_1, P_2)$ be a bihamiltonian structure, and denote
\[D_1=D_{P_1},\quad D_2=D_{P_2},\quad d_1=d_{P_1},\quad d_2=d_{P_2}.\]
These operators satisfy the identities
\[D_1 D_2+D_2 D_1=0,\quad d_1 d_2+d_2 d_1=0,\]
so we obtain two double complexes $(\AA, D_1, D_2)$ and $(\FF, d_1, d_2)$.

\begin{dfn}\label{ky-15}
Let $(C, \p_1, \p_2)$ be the double complex $(\AA, D_1, D_2)$ or $(\FF, d_1, d_2)$, its bihamiltonian cohomology is defined as
\begin{equation}
BH^p_d(C, \p_1, \p_2)=\frac{C^p_d\cap\Ker\p_1\cap \Ker \p_2}{C^p_d\cap\Img(\p_1\,\p_2)}, \quad p, d\ge 0.
\end{equation}
We often denote it by $BH^p_d(C)$ for short.
\end{dfn}

\begin{rmk}
Due to Lemma \ref{lem-36} and Theorem \ref{zh-thm-36}, the above definition of bihamiltonian cohomologies is in agreement with the one given in \eqref{zh-37}, \eqref{zh-37-b} for finite dimensional bihamiltonian structures
when $d\ge 2$.
\end{rmk}

\begin{lem}\label{lem-43}
Let $(C, \p_1, \p_2)$ be the double complex $(\AA, D_1, D_2)$ or $(\FF, d_1, d_2)$. Define
\[C_\lambda=C\otimes \R[\lambda],\quad \p_\lambda=\p_2-\lambda\,\p_1,\]
then $(C_\lambda, \p_\lambda)$ is a complex. Moreover, if $P_1, P_2$ is a Hamiltonian structure of hydrodynamic type, then for any $p\in\N$ and $d\ge2$, we have
\begin{equation}
BH^p_d(C, \p_1, \p_2)\cong H^p_d(C_\lambda, \p_\lambda). \label{eq-iso}
\end{equation}
\end{lem}

\begin{prf}
We embed $C$ into $C_\lambda$ in the natural way. It's easy to see that if $z\in C^p_d\cap\Ker\p_1 \cap \Ker\p_2$,
then $z\in (C_\lambda)^p_d\cap\Ker \p_\lambda$,
so there is an embedding
\[j: C^p_d\cap\Ker\p_1 \cap \Ker\p_2 \to (C_\lambda)^p_d\cap\Ker \p_\lambda.\]
We are to show that this embedding induces the isomorphism \eqref{eq-iso}.

Let $\pi$ be the projection 
\[\pi:  (C_\lambda)^p_d\cap\Ker \p_\lambda \to H^p_d(C_\lambda, \p_\lambda),\] 
we first show
that the composition $\pi\, j$ is surjective. Let $x\in(C_\lambda)^p_d\cap\Ker\p_\lambda$, we need to prove that there exists
$y\in (C_\lambda)^{p-1}_{d-1}$ such that $x-\p_\lambda y\in \Img j$.
Suppose $x$ takes the following form
\[x=x_0+x_1\,\lambda+\cdots+x_m\,\lambda^m,\quad\mbox{where } x_k\in C^p_d.\]
The condition $\p_\lambda x=0$ implies that
\begin{equation}\label{zh-38}
0=\p_2\, x_0,\quad \p_1\, x_0=\p_2\, x_1, \quad \cdots, \quad \p_1\, x_{m-1}=\p_2\, x_m,\quad \p_1\, x_m=0.
\end{equation}
Since $H^p_{d}(C, \p_1)\cong0$, we know that there exists $y_{m-1}\in C^{p-1}_{d-1}$ such that
\[x_m=-\p_1\, y_{m-1},\]
then we have $\p_1\,x_{m-1}=\p_1\,\p_2\, y_{m-1}$, so there exists $y_{m-2}\in C^{p-1}_{d-1}$ such that
\[x_{m-1}=\p_2\, y_{m-1}-\p_1\, y_{m-2}.\]
By induction, one can show the existence of  $y_{k-1}\in C^{p-1}_{d-1}$ such that
\begin{equation}\label{zh-39}
x_{k}=\p_2\, y_{k}-\p_1\,y_{k-1},\quad  k=1, \dots, m-1.
\end{equation}
Now let
\[y=y_0+y_1\,\lambda+\cdots+y_{m-1}\,\lambda^{m-1}\in (C_\lambda)^{p-1}_{d-1},\]
then from \eqref{zh-38}, \eqref{zh-39} it follows that 
\[x-\p_\lambda\,y=x_0-\p_2\, y_0\in \Img j.\]
The surjectivity is proved.

Next we need to show that $\Ker(\pi\, j)=C^p_d\cap\Img (\p_1\,\p_2)$. If $z=\p_1\,\p_2\,w$, then
\[z=\p_\lambda(-\p_1\,w),\]
so $C^p_d\cap\Img (\p_1\,\p_2)\subset \Ker(\pi\, j)$. Conversely, suppose $z\in C^p_d\cap\Ker \p_1 \cap \Ker \p_2$ is exact in
$H^p_d(C_\lambda, \p_\lambda)$, i.e. there exists $y\in (C_\lambda)^{p-1}_{d-1}$ such that $z=\p_\lambda\,y$.  By using the fact that $H^{p-1}_{d-1}(C,\p_1)\cong0$ and a similar argument as we used above in the proof of the surjectivity of $\pi\, j$, it is easy to show that we can choose this $y$ such that it belongs to  $C$. Then we have
\[z=\p_2\,y,\quad 0=\p_1\,y.\]
Since $H^{p-1}_{d-1}(C,\p_1)\cong 0$, there exists $w\in C^{p-2}_{d-2}$ such that $y=-\p_1\,w$. From this it follows that  $z=\p_1\,\p_2\,w$,  so
we have $\Ker(\pi\, j)\subset \Img (\p_1\,\p_2)$.

Finally,
\[H^p_d(C_\lambda, \p_\lambda)=\Img (\pi\,j)\cong \frac{C^p_d\cap\Ker\p_1 \cap \Ker \p_2}{\Ker (\pi\,j)}=BH^p_d(C,\p_1, \p_2).\]
The lemma is proved.
\end{prf}

\begin{rmk}
A similar description of $BH^p(C, \p_1, \p_2)$ was given in \cite{B}, in which the role of the complex
$(C_\lambda, \p_{\lambda})$ is replaced by a bicomplex $C^{\bullet, \bullet}$.
\end{rmk}

\begin{cor}\label{cor-44}
If $P_1, P_2$ are Hamiltonian structures of hydrodynamic type, then for any $p\in\N$ and $d\ge2$ there exists
a long exact sequence:
\begin{equation}
\cdots\xrightarrow{\ \ \ }BH^p_d(\AA)\xrightarrow{\ \ \ }
BH^p_d(\FF)\xrightarrow{\ \ \ }BH^{p+1}_d(\AA)\xrightarrow{\ \ \ }BH^{p+1}_{d+1}(\AA)\xrightarrow{\ \ \ }\cdots.\label{bhc-long}
\end{equation}
\end{cor}
\begin{prf}
Denote
\begin{align*}
\AA_\lambda&=\AA\otimes \R[\lambda],\quad  D_\lambda=D_2-\lambda\,D_1,\\
\FF_\lambda&=\FF\otimes \R[\lambda],\quad  d_\lambda=d_2-\lambda\,d_1,
\end{align*}
then we have a short exact sequence
\[0\xrightarrow{\ \ \ }(\AA_\lambda/\R[\lambda], D_\lambda)
\xrightarrow{\ \p\ }(\AA_\lambda, D_\lambda)
\xrightarrow{\ \int\ }(\FF_\lambda, d_\lambda)\xrightarrow{\ \ \ }0.\]
This short exact sequence implies a long exact sequence
\begin{equation}
\cdots\xrightarrow{\ \ \ }H^p_d(\AA_\lambda)\xrightarrow{\ \ \ }
H^p_d(\FF_\lambda)\xrightarrow{\ \ \ }H^{p+1}_d(\AA_\lambda/\R[\lambda])
\xrightarrow{\ \ \ }H^{p+1}_{d+1}(\AA_\lambda)\xrightarrow{\ \ \ }\cdots, \label{long2}
\end{equation}
which is isomorphic to \eqref{bhc-long} when  $d\ge 2$ due to Lemma \ref{lem-43}. The corollary is proved.
\end{prf}

Let $(P_1, P_2)$ be a semisimple bihamiltonian structure of hydrodynamic type. A deformation of $(P_1, P_2)$
is a pair $(\tilde{P}_1, \tilde{P}_2)\in \FF^2\times\FF^2$ which satisfies
\begin{align*}
&\nu(\tilde{P}_1-P_1)>\nu(P_1),\quad \nu(\tilde{P}_2-P_2)>\nu(P_2),\\
&[\tilde{P}_1, \tilde{P}_1]=[\tilde{P}_1, \tilde{P}_2]=[\tilde{P}_2, \tilde{P}_2]=0.
\end{align*}
Here the valuation $\nu$ is defined in \eqref{ky-9}. Note that, according to Theorem \ref{zh-thm-36}, there always exists a Miura type transformation that transforms
$\tilde{P}_1$ to $P_1$. Thus without loss of generality we only need to consider deformations of the form
$(P_1,\tilde{P}_2)$, hence the
deformation theory of the bihamiltonian structure $(P_1,  P_2)$ is controlled by the DGLA $(\mathrm{Ker}(d_1), [\ \,,\ ], d_2)$ whose cohomology groups
are just the bihamiltonian cohomologies. We call a pair $(Q_1, Q_2)\in \FF^2\times\FF^2$ an infinitesimal deformation of the bihamiltonian structure $(P_1, P_2)$ if $\nu(Q_1)>\nu(P_1),\ \nu(Q_2)>\nu(P_2)$ and 
\begin{equation*}
[P_1, Q_1]=[P_2, Q_2]=[P_1, Q_2]+[P_2, Q_1]=0.
\end{equation*}  
Also due to Theorem \ref{zh-thm-36} we only need to consider infinitesimal deformations of the form $(0, Q)$.

It is easy to see that
\begin{itemize}
\item[i)] The group $BH^2_{\ge 2}(\FF)$ characterizes  the equivalence classes of infinitesimal deformations of the bihamiltonian structure $(P_1, P_2)$.
\item[ii)] If $BH^3_{\ge 2\delta}(\FF)$ vanishes, then any infinitesimal deformation of $(P_1, P_2)$ can be extended to a deformation of it. Here $\delta$ is the lowest degree of classes in $BH^2_{\ge 2}(\FF)$.
\end{itemize}

In our previous works \cite{LZ1, DLZ1}, we prove that:
\begin{equation}\label{ky-10}
BH^2_d(\FF)\cong 0, \quad\mbox{if $d\ge2$ and $d\ne 3$},
\end{equation}
and
\begin{equation}
BH^2_3(\FF)\cong \bigoplus_{i=1}^n C^\infty(\R).
\end{equation}
This knowledge of the second bihamiltonian cohomologies  enables us to classify the infinitesimal deformations of $(P_1, P_2)$.
In particular, for a bihamiltonian structure of hydrodynamic type \eqref{ky-5} represented via the loop space 
approach, as we explained in Remark \ref{ky-8}, we only need to consider its infinitesimal deformations given by  \eqref{ky-6}, \eqref{ky-6-b}. In order to study 
the existence of a deformation of $(P_1, P_2)$ with a given infinitesimal one, we need to compute the third cohomologies $BH^3_{\ge6}(\FF)$ and to show its triviality (see \cite{LZ1} for more details).

The computation carried out in \cite{LZ1, DLZ1} to prove the above result on the second bihamiltonian cohomologies
is quite involved. This is mainly due to the fact that the definition of the differentials $d_1, d_2$ on $\FF$ are complicated,
so it is difficult to perform similar computations  to study the third bihamiltonian cohomologies. Fortunately, the long exact sequence \eqref{bhc-long} given in Corollary \ref{cor-44} helps us to overcome this difficulty.  From this long exact sequence it follows that the triviality of $BH^3_{d\ge 6}(\FF)$ can be
proved by showing that $BH^3_{d\ge 6}(\AA)\cong0$ and  $BH^{4}_{d\ge 6}(\AA)\cong0$. Note that the definition of the
differentials $D_1$, $D_2$ on $\AA$ are much simpler than that of the differentials $d_1, d_2$
on $\FF$, so it is easier to compute the bihamiltonian cohomologies $BH^p_d(\AA)$ than to compute  $BH^p_d(\FF)$.

In the following section, we illustrate the above mentioned method of computing bihamiltonian cohomologies by proving Theorem \ref{thm-zh1} and Theorem \ref{thm-zh2} for the  bihamiltonian structure \eqref{zh-19}, \eqref{zh-20} of the dispersionless KdV hierarchy. 

\section{The third bihamiltonian cohomology for the dispersionless KdV hierarchy}
\label{sec-5}

The following subsections are devoted to  prove Theorem \ref{thm-zh1} and Theorem \ref{thm-zh2}.
\subsection{Some preparations}
For the bihamiltonian structure of the dispersionless KdV hierarchy we have $M=\R$, so we can omit the index $i$ in the notations such as $u^{i,s}, \theta_i^s$. The bihamiltonian structure
\eqref{zh-19}, \eqref{zh-20} can be represented by the bivectors
\begin{equation}\label{bh-dkdv}
P_1=\frac12\int \theta\,\theta^1,\quad P_2=\frac12\int u\,\theta\,\theta^1,
\end{equation}
so we have
\begin{align*}
D_1&=\sum_{s\ge0}\theta^{s+1}\frac{\p}{\p u^s},\\
D_2&=\sum_{s\ge0}\left(\p^s\left(u\,\theta^1+\frac12\,u^1\,\theta\right)\frac{\p}{\p u^s}
+\p^s\left(\frac12\,\theta\,\theta^1\right)\frac{\p}{\p \theta^s}\right).
\end{align*}
The following identities are useful and easy to prove, so we list them here and omit the proof.
\begin{lem}\label{lem-51}
\begin{align*}
[\frac{\p}{\p u^k},\ D_1]=&0,\\
[\frac{\p}{\p \theta^k},\ D_1]=&\frac{\p}{\p u^{k-1}},\\
[\frac{\p}{\p u^k},\ D_2]=&
\sum_{l\ge k-1}\left(\binom{l}{k}+\frac12\binom{l}{k-1}\right)\theta^{l+1-k}\frac{\p}{\p u^l},\\
[\frac{\p}{\p \theta^k},\ D_2]=&
\sum_{l\ge k-1}\left(\binom{l}{k-1}+\frac12\binom{l}{k}\right)u^{l+1-k}\frac{\p}{\p u^l}\\
&+\frac12\sum_{l\ge k-1}\left(\binom{l}{k}-\binom{l}{k-1}\right)\theta^{l+1-k}\frac{\p}{\p \theta^l}.
\end{align*}
\end{lem}

Suppose $Q\in\AA^p_d\ (d\ge2)$ satisfies
\[D_1 Q=0,\quad D_2 Q=0.\]
To prove the triviality of $BH^p_d(\AA)$, we need to find $X\in\AA^{p-2}_{d-2}$ such that
\[Q=D_1 D_2 X.\]
According to Lemma \ref{lem-36}, there exists $F\in \AA^{p-1}_{d-1}$, such that
\[Q=D_1 F,\quad  D_1 D_2 F=0.\]
So in order to prove Theorem \ref{thm-zh1} we need to do the following:

\textit{For any given $(p,d)=(3,\ge6) \mbox{ or } (4, \ge6)$ and $F\in\AA^{p-1}_{d-1}$ satisfying $D_1 D_2 F=0$,
find $X\in\AA^{p-2}_{d-2}$ such that $D_1 F=D_1 D_2 X$.}

Our proof is in the spirit of spectral sequence, especially the one of a filtered complex.
We first introduce a filtration on $\AA$
\[0\subset\AA^{(0)}\subset\AA^{(1)}\subset\AA^{(2)}\subset\cdots\subset\AA,\]
where $\AA^{(k)}$ is the differential polynomial algebra on $J^k(\hM)$.
Suppose
\[Q\in \Ker D_1\cap \Ker D_2 \cap \AA^{p,(N+1)}_d,\]
if we can find $X\in\AA^{p-2,(N-1)}_{d-2}$ such that
\[\tilde{Q}=Q-D_1 D_2 X\in\AA^{p,(N)}_d,\]
then $\tilde{Q}$ gives a new representative of the cohomology class of $Q$, which depends on
less variables. By induction on $N$, we will show that one can choose a representative $Q$ with smallest $N$
(in the  $BH^3$ case, the smallest $N$ is $1$, while in the $BH^4$ case it is $2$). For this $Q$, the equations
$D_1 Q=0$ and $D_2 Q=0$ are easy to solve, and in this way one is able to
compute the  bihamiltonian cohomologies. So the most important step in our computation of the bihamiltonian cohomologies is the following:

\textit{For any given $(p,d)=(3,\ge6) \mbox{ or } (4, \ge6)$ and $F\in\AA^{p-1,(N)}_{d-1}$ satisfying $D_1 D_2 F=0$, find $X\in\AA^{p-2,(N-1)}_{d-2}$ such that $D_1 F-D_1 D_2 X\in\AA^{p,(N)}_d$.}

In this subsection, we prove some general results which will be used in the next two subsections.

\begin{lem}\label{lem-52}
If $F\in\AA^{(N)}$ with  $N\ge 1$ satisfies
\[\frac{\p F}{\p u^N}=0,\]
then there exists $X\in\AA^{(N-1)}$ such that $F-D_1 X\in\AA^{(N-1)}$.
\end{lem}
\begin{prf}
Let us represent $F$ in the form
\[F=\theta^NF_N+F_0,\]
where $F_N, F_0\in\AA^{(N-1)}$.  Then one can check that
\[X=\p_{u^{N-1}}^{-1}F_N\]
fulfills the condition $F-D_1 X\in\AA^{(N-1)}$. Here and in what follows we use  
the notation
\[\p_{t}^{-1}A(t)=\int_0^{t} A(\tau)\,\mathrm{d}\tau.\]
Note that $F_N$ is always a formal power series of $\epsilon$ whose coefficients are smooth function
or polynomials of $u^{N-1}$, so the operator
$\p_{u^{N-1}}^{-1}$ is well defined. 
\end{prf}

\begin{lem}\label{lem-53}
If $F\in\AA^{(N)}$ with $N\ge 1$ satisfies the condition
\[D_1 D_2 F\in\AA^{(N+1)},\]
then there exists $X\in\AA^{(N-1)}$ such that $\tilde{F}=F-D_1 X$
 has the form $\tilde{F}=\theta\,A+B$, where $A\in\AA^{(N)}$, $B\in\AA^{(N-1)}$.
\end{lem}
\begin{prf}
From the given assumptions on $F$ it follows that
\[\frac{\p}{\p u^{N+2}}\left(D_1 D_2 F\right)=0,\quad  \frac{\p}{\p \theta^{N+2}}\left(D_1 D_2 F\right)=
\frac12\theta\frac{\p F}{\p u^N}.\]
We write $F$ as $\theta\,A+B_1$, where $B_1$ is independent of $\theta$, then the above equations
and the condition $D_1 D_2 F\in\AA^{(N+1)}$ imply that
\[\frac{\p B_1}{\p u^N}=0.\]
Then it follows from Lemma \ref{lem-52} that there exists $X\in\AA^{(N-1)}$ such that $B=B_1-D_1 X\in\AA^{(N-1)}$.
The lemma is proved.
\end{prf}

\begin{lem}\label{lem-54}
Let $R$ be an $\R$-algebra, $g\in R[x_1, \dots, x_\ell]$, then for any $a_0, \dots, a_\ell\in \R_{>0}$,
the equation
\begin{equation}
a_0 f+a_1 x_1 \frac{\p f}{\p x_1}+\cdots+a_\ell x_\ell \frac{\p f}{\p x_\ell}=g \label{eq-pde}
\end{equation}
has a unique solution $f\in R[x_1, \dots, x_\ell]$.
\end{lem}
\begin{prf}
Let us first prove the existence of $f$. Every element of $R[x_1, \dots, x_\ell]$ can be written as a finite sum of
monomials, so we only need to prove the case when $g$ is a monomial. 
Suppose
\[g=r\,x_1^{m_1}\cdots x_\ell^{m_\ell},\quad r\in R, \quad m_1, \dots, m_\ell \in \N,\]
then the equation \eqref{eq-pde} has a solution
\[f=\frac{r}{a_0+m_1\,a_1+\cdots+m_\ell\,a_\ell}\,x_1^{m_1}\cdots x_\ell^{m_\ell}.\]
Here the condition that $a_k>0$ ensures that the denominator does not vanish.

To prove the uniqueness, we only need to show that the homogeneous equation of \eqref{eq-pde}
has only zero solution. Suppose $f\in R[x_1, \dots, x_\ell]$ is a nonzero solution to \eqref{eq-pde}
with $g=0$, then we can write $f$ as
\[f=\sum_{m_1,\dots,m_\ell\ge0} r_{m_1,\dots,m_\ell}x_1^{m_1}\cdots x_\ell^{m_\ell}.\]
Let $K$ be the degree of $f$, then we know that there exist $m_1, \dots, m_\ell$ such that
\[m_1+\cdots+m_\ell=K,\quad  r_{m_1,\dots,m_\ell}\ne 0.\]
We denote $\underline{m}=(m_1, \dots, m_\ell)$, and
\[\p^{\underline{m}}=\p_{x_1}^{m_1}\cdots\p_{x_\ell}^{m_\ell},\quad  x^{\underline{m}}=x_1^{m_1}\cdots x_\ell^{m_\ell},\]
then compute the action of $\p^{\underline{m}}$ on the equation \eqref{eq-pde}.
From the choice of $\underline{m}$ it follows that
\[(a_0+m_1\,a_1+\cdots+m_\ell\,a_\ell)r_{m_1,\dots,m_\ell}\p^{\underline{m}}(x^{\underline{m}})=0.\]
This is impossible, since every term on the left hand side does not vanish, so the lemma is proved.
\end{prf}

\begin{rmk}\label{rmk-55}
By using the same method, one can prove that if $g$ has no constant term, i.e. $g(0, \dots, 0)=0$,
then the following equation
\[a_1 x_1 \frac{\p f}{\p x_1}+\cdots+a_\ell x_\ell \frac{\p f}{\p x_\ell}=g\]
has a solution when $a_1, \dots, a_\ell>0$. This fact is also used in the following sections.
But in this case, the solution is not unique, since arbitrary constant is a solution to the corresponding homogeneous equation.
\end{rmk}

\begin{lem}\label{lem-55}
Suppose $F=\theta\,A+B\in\AA^{(N)}$ with $N\ge 1$, where
\[\frac{\p A}{\p u^N}\in\AA^{(N-1)},\quad  \frac{\p A}{\p \theta^N}=0, \quad  B\in\AA^{(N-1)},\]
then there exists $X, Y\in\AA^{(N-1)}$ such that $F-D_2 X+D_1 Y\in\AA^{(N-1)}$.
\end{lem}
\begin{prf}
The proof is similar to the one for Lemma \ref{lem-52}, so we only present 
\[X=\p_{u^{N-1}}^{-1}\left(2\frac{\p A}{\p u^N}\right),
\quad  Y=\p_{u^{N-1}}^{-1}\left(u\,\frac{\p X}{\p u^{N-1}}-\frac12\theta\,\frac{\p X}{\p \theta^{N-1}}\right),\]
and omit the details. 
\end{prf}

Now we give  the proof of triviality of $BH^2(\AA)$ to illustrate the usage of these lemmas.
\begin{thm}\label{zh-27-5}
For the second bihamiltonian cohomologies we have
\[BH^2_d(\AA)=\left\{\begin{array}{ll}
0, &\quad  d=0,\\
\{a(u)\theta\theta'\,|\,a(u)\in C^{\infty}({\mathbb{R})}\},  &\quad  d=1,\\
0, &\quad d\ge2.
\end{array}\right.\]
\end{thm}
\begin{prf}
In the $d=0$ or $1$ case, since there is no $d-2$ degree component in $\AA$, so $BH^2_d(\AA)$ is just
$\Ker D_1 \cap \Ker D_2 \cap \AA^2_d$, which is easy to compute. So we only give the proof for the
$d\ge 2$ cases.

Suppose $F\in\AA^{1,(N)}_d$ satisfies $D_1 D_2 F=0$. Due to Lemma \ref{lem-53} we can assume that
$F$ has the following form:
\[F=a\,\theta+b\]
where $a\in \AA^{0,(N)}=\A^{(N)}$, $b\in \AA^{1,(N-1)}$.  A simple calculation shows that 
\[\frac{\p^2}{\p \theta^N \p \theta^{N+1}}\left(D_1 D_2 F\right)=\left(N+\frac12\right)
u^1\theta  \frac{\p ^2 a}{\p u^N\p u^N},\]
so from the vanishing of $D_1 D_2 F$ it follows that $a$ must take the  form
\[a=a_0(u, \dots, u^{N-1})u^{N}+a_1(u, \dots, u^{N-1}),\]
thus $F$ satisfies the conditions of Lemma \ref{lem-55}.
By using this lemma, we can find $X, Y\in\AA^{(N-1)}$ such that
$\tilde{F}=F-D_2 X+D_1 Y\in\AA^{1,(N-1)}$.

By induction on $N$, we can reduce $F$ to an element of the space $\AA^{1,(0)}$. But $\AA^{1,(0)}=\AA^1_0$, so
$Q=D_1 F\in\AA^2_1$ cannot be an element of $\AA^2_{d\ge 2}$ unless it vanishes.
The theorem is proved.
\end{prf}

\subsection{The bihamiltonian cohomology $BH^3(\AA)$}

We are to prove the following theorem in this subsection.
\begin{thm}\label{thm-6}
The bihamiltonian cohomology $BH^3_d(\AA)$ is given by 
\[BH^3_d(\AA)=\left\{\begin{array}{ll}
0, &\quad d=0, 1,2,\\
\{a(u)\theta\theta^1\theta^2\,|\,a(u)\in C^{\infty}(\mathbb{R})\}, &\quad d=3,\\
0,  &\quad  d\ge 4.
\end{array}\right.\]
\end{thm}

\begin{prf}
For any given $Q\in \Ker D_1\cap \Ker D_2 \cap \AA^{3,(N+1)}, \ N\ge 1$, we need to prove the existence of 
$X\in \AA^1$ such that $Q-D_1 D_2 X$ has the expression $a(u)\theta\theta^1\theta^2$. To this end, let us first
find $F\in\AA^{2,(N)}$ such that $Q=D_1 F$, then $F$ satisfies $D_1 D_2 F=0$.
We are to show that, by modifying $F$ in the way $F\mapsto F-D_2 X+D_1 Y$ with certain $X, Y\in\AA^1$, 
we can reduce $F$ to an element of the space $\AA^{2,(1)}$.

Denote $Z=D_1 D_2 F$, then by using the fact that $Z=0$
and the result of  Lemma \ref{lem-53}  we can reduce $F$ to the form
\begin{equation}\label{zh-2}
F=\sum_{s=1}^N A_s\theta\theta^s+B,
\end{equation}
where $A_s\in\A^{(N)},\ B\in\AA^{2,(N-1)}$.

Let us assume $N\ge 2$ and denote the coefficients of $\theta\theta^{N-s}\theta^N\theta^{N+1}$
in $Z$ by
\begin{equation}\label{zh-3}
Z_s=\frac{\p^4 Z}{\p\theta^{N+1}\p\theta^N\p\theta^{N-s}\p\theta},\quad   s=1,2,\dots,N-1,
\end{equation}
then by using Lemma \ref{lem-51}, one can easily obtain that
\[Z_1\equiv \frac{N^2}{2} u^2\,\frac{\p^2 A_N}{\p u^N \p u^N}  \  (\textrm{mod\ } u^1).\]
Since $Z=0$, we arrive at
\[u^1 \left|\frac{\p^2 A_N}{\p u^N \p u^N}\right..\]
By using this fact, one can further obtain that
\[0=Z_2\equiv  \frac{(N-1)^2}{2} u^2\,\frac{\p^2 A_N}{\p u^{N-1} \p u^N}  \  (\textrm{mod\ } u^1 ),\]
so we have
\[u^1 \left|\frac{\p^2 A_N}{\p u^{N-1} \p u^N}\right..\]
By induction on $s$, one can prove that
\[u^1 \left|\frac{\p^2 A_N}{\p u^{N-s+1} \p u^N}\right.,\quad  s=1, \dots, N-2.\]
When $s=N-1$, we have
\[0=Z_{N-1}\equiv  2 u^2\,\frac{\p^2 A_N}{\p u^{2} \p u^N}+(N-1) \frac{\p A_N}{\p u^N}  \  (\textrm{mod\ } u^1 ).\]
Denote $f=\frac{\p A_N}{\p u^N}$, then
\[2 u^2\frac{\p f}{\p u^2}+(N-1)f\equiv 0\  (\textrm{mod\ } u^1),\]
By using Lemma \ref{lem-54} and the above equation we know that
\[u^1 \left|\frac{\p A_N}{\p u^N}\right.,\]
so the differential polynomial $A_N$ can be put into the form
\[A_N=u^1 a+b, \quad \mbox{with }  a\in\A^{(N)},\quad  b\in\A^{(N-1)}.\]
Then Lemma \ref{lem-52} enables us to modify the term 
$b\,\theta\,\theta^N$ of $F$ by adding certain  $D_1 X\in\AA^{2,(N)}$ 
such that the resulting bivector belongs to $\AA^{2,(N-1)}$, and we can combine it into the term $B$
in the expression \eqref{zh-2} of $F$, so we can assume that 
\[A_N=u^1 a.\]
Let us choose
\[X=-(N+\frac12)^{-1} \p^{-1}_{u^N}\left(a\right)\theta,\quad   Y=u X,\]
then it is easy to see that after the modification $F\mapsto F-D_2 X+D_1 Y $
we can reduce $F$ to the form \eqref{zh-2} with $A_N=0$ and $B\in\AA^{2,(N-1)}$.
Now the functions $Z_s$ defined in \eqref{zh-3} have the expression
\[0=Z_s=-\left(N+\frac12\right) u^1 \frac{\p^2 A_{N-s}}{\p u^N\p u^N},\quad  s=1,\dots,N-1,\]
so the vector $A=\sum\limits_{s=1}^{N-1} A_s(u,u^1,\dots,u^N) \theta^s$ satisfies the condition
\[\frac{\p A}{\p u^N}\in\AA^{(N-1)}, \quad   \frac{\p A}{\p \theta^N}=0,\]
and it follows from Lemma \ref{lem-55} that we can find $X, Y\in\AA^{1,(N-1)}$ such that
\[F-D_2 X+D_1 Y\in\AA^{2,(N-1)}.\]

By induction on $N$, we can prove the existence of $X, Y\in\AA^1$ such that 
$F-D_2 X-D_1 Y\in\AA^{2,(1)}$, so we can reduce $F$ to the form
\[F=A_1(u,u^1) \theta\theta^1.\]
Then the trivector $D_1 F$ has the expression 
$D_1 F=\frac{\p A_1(u,u^1)}{\p u^1} \theta\theta^1\theta^2$.
Denote $\frac{\p A_1(u,u^1)}{\p u^1}=u^1 b(u,u^1)+a(u)$, and define $X=-\frac23 \p^{-2}_{u^1}\left(b(u,u^1)\right)\,\theta$, then we have
\[Q\sim D_1 F-D_1 D_2 X=a(u) \theta\theta^1\theta^2.\]
Thus we proved the theorem.
\end{prf}

\subsection{The bihamiltonian cohomology $BH^4(\AA)$} 

In this subsection, we continue to compute the bihamiltonian cohomology $BH^4(\AA)$.

\begin{thm}\label{thm-zh7}
The bihamiltonian cohomology $BH^4(\AA)$ is trivial. \end{thm}

\begin{prf}
Given $Q\in \Ker D_1\cap \Ker D_2 \cap \AA^{4,(N+1)},\ N\ge 2$, we need to prove the existence of 
$X\in \AA^2$ such that $Q=D_1 D_2 X$. As 
we did in the computation of the bihamiltonian cohomology $BH^3(\AA)$, we can
find $F\in\AA^{3,(N)}$ such that $Q=D_1 F$, and $F$ satisfies $D_1 D_2 F=0$.
Let us  prove the theorem in three steps. In the first step, we prove that when $N\ge 4$ we can reduce $F\in\AA^{3,(N)}$, by modifying it in the way $F\mapsto F-D_2 X+D_1 Y$ with certain $X, Y\in\AA^2$,
to an element of the space $\AA^{3,(3)}$. In the second step, we prove that an element $F$ of $\AA^{3,(3)}$ satisfying $D_1 D_2 F=0$ can be reduced to an element of the space  $F\in\AA^{3,(2)}$. In the last step, we prove that 
for any element $F$ of $\AA^{3,(2)}$ satisfying $D_1 D_2 F=0$, one can always represent $D_1 F$ as $D_1 D_2 X$ for a certain $X\in\AA^2$. 

\vskip 0.2cm
\noindent Step 1: For the case when $N\ge 4$.

Denote $Z=D_1 D_2 F$, then by using Lemma \ref{lem-53} and $Z=0$ we know that $F$ can
be assumed to have the form
\begin{equation}\label{zh-5}
F=\sum_{1\le p<q\le N} A_{p,q}\theta\theta^p\theta^q+B,
\end{equation}
where $A_{p,q}\in\A^{(N)}$, $B\in\AA^{3,(N-1)}$. We define
\begin{equation}\label{zh-6}
Z_s=\frac{\p^4 Z}{\p\theta^{N+1}\p\theta^N\p\theta^{N-s}\p\theta},\quad  s=1,2,\dots,N-1,
\end{equation}
then by using Lemma \ref{lem-51} and induction on $p$ we have
\begin{align*}
& \frac{\p Z_1}{\p\theta^{p}}\equiv \frac{(p+1)(p+2)(2\,p+3)}{12} u^3\,\frac{\p^2 A_{N-1,N}}{\p u^{p+2} \p u^N}  \ 
(\textrm{mod\ } u^1, u^2),\quad 2\le p\le N-2,\\
&  \frac{\p Z_1}{\p\theta^{1}}\equiv \frac52 u^3 \frac{\p^2 A_{N-1,N}}{\p u^3\p u^N}+\frac{3 N-4}2\, \frac{\p A_{N-1,N}}{\p u^N} 
\  (\textrm{mod\ } u^1, u^2).
\end{align*}
So we have $\frac{\p A_{N-1,N}}{\p u^N}\equiv 0 \ (\textrm{mod\ } u^1, u^2)$. A similar argument as we gave in the proof of Theorem \ref{thm-6} shows that $A_{N-1,N}$ can be represented as
\[A_{N-1,N}=u^1 B_1+u^2 B_2,\quad B_1, B_2 \in \A^{(N)}.\]

Now we assume $X, Y \in\AA^2$ have the form
\[X=x_0\theta\theta^N+x_1\theta\theta^{N-1},\quad  Y=u X,\quad  x_0, x_1\in\A^{(N)},\]
and we consider the following modification of $F$ 
\[\tilde{F}=F-D_2 X+D_1 Y.\]
Then $\tilde{F}\in\AA^{3,(N)}$ and the coefficient of $\theta\theta^{N-1}\theta^N$ in  $\tilde{F}$ has the expression
\begin{align*}
&\frac{\p^3\tilde{F}}{\p\theta\p\theta^{N-1}\p\theta^N}\\
&=u^1 \left(\frac{2 N+1}2 \frac{\p x_1}{\p u^N}-\frac{2N-1}2 \frac{\p x_0}{\p u^{N-1}}-B_1\right) -u^2 \left(\frac{N^2}2 \frac{\p x_0}{\p u^N}+B_2\right).
\end{align*}
From the above formula it follows the existence of differential polynomials $x_0, x_1$ such that $\frac{\p^3\tilde{F}}{\p\theta\p\theta^{N-1}\p\theta^N}=0$ and we have reduced 
$F$ to the form \eqref{zh-5} with $A_{N-1,N}=0$. 

We are to prove by induction that $F$ can actually be reduced to the form \eqref{zh-5} with $A_{p,N}=0, \
p=1,2,\dots,N-1$. Assume that $F$ has been reduced to the form
\begin{equation}\label{zh-7}
F=\sum_{1\le p<q\le N-1} A_{p,q}\theta\theta^p\theta^q+\sum_{1\le p\le m} A_{p,N}\theta\theta^p\theta^N+B
\end{equation}
with $A_{p,q}\in\A^{(N)}$ and $B\in\AA^{3,(N-1)}$ for a certain $ 1\le m\le N-2$. Then the vectors $Z_s$ defined in \eqref{zh-6} have the following properties
(which can be proved by induction on $s$)
\begin{align*}
0=&\frac{\p Z_s}{\p \theta^p}\equiv -\frac{(N-s+1)^2}2 u^2 \frac{\p^2 A_{p,N}}{\p u^{N-s+1}\p u^N}\  (\textrm{mod\ } u^1),\\
& \qquad\qquad\qquad\qquad   s=1\dots, N-m-1,\quad p=1,\dots, m.
\end{align*}
These relations imply that
\begin{equation}\label{divide}
u_1 \left|\frac{\p^2 A_{p,N}}{\p u^{N-s+1}\p u^N}\right.,\quad s=1\dots, N-m-1,\quad p=1,\dots, m.
\end{equation}
Let us consider the cases $m\ge 2$ and $m=1$ separately.

\noindent\ i)\ The $m\ge 2$ case. 
From \eqref{divide} it follows that (here an induction on $p$ is needed)
\begin{align*}
& \frac{\p Z_{N-m}}{\p\theta^{p}}\equiv \frac{(p+1)(p+2)(2\,p+3)}{12} u^3\,\frac{\p^2 A_{m,N}}{\p u^{p+2} \p u^N}  \ (\textrm{mod\ } u^1, u^2),
\quad  2\le p\le m-1,\\
&  \frac{\p Z_{N-m}}{\p\theta^{1}}\equiv \frac52 u^3 \frac{\p^2 A_{m,N}}{\p u^3\p u^N}+\frac{2 N+m-3}2\, \frac{\p A_{m,N}}{\p u^N}  \ (\textrm{mod\ } u^1, u^2).
\end{align*}
The last equation together with Lemma \ref{lem-54} imply that 
\[\frac{\p A_{m,N}}{\p u^N}\equiv 0 \ (\textrm{mod\ } u^1, u^2).\]
So we can represent the differential polynomial  $A_{m,N}$  in the form
\begin{equation}
A_{m,N}=u^1 B_1(u,u^1,\dots,u^N)+u^2 B_2(u,u^2,\dots,u^N)
\end{equation}
after eliminating a remainder term that does not depend on $u^N$ by using Lemma \ref{lem-52}.
Note that the differential polynomial $B_2$ is chosen to be independent of $u^1$.
By using \eqref{divide} with $p=m$ we know that 
\begin{equation}\label{zh-8}
\frac{\p^2 B_2}{\p u^{m+2}\p u^N}= \dots=\frac{\p^2 B_2}{\p u^{N}\p u^N}=0.
\end{equation}

Now let us try to find $X, Y\in\AA^{2}$ of the form
\begin{align}
X=&x_0\theta\theta^N+x_1\theta\theta^m+x_2\theta\theta^{m+1},\quad  x_0, x_1, x_2\in\A^{(N)},\label{zh-11-27-1}\\
Y=&u X-\sum_{k=m}^{N-1}y_k\theta\theta^k,\quad  y_k\in\A^{(N-1)},\label{zh-11-27-2}
\end{align}
such that the following modification of $F$
\begin{equation}
\tilde{F}=F-D_2 X+D_1 Y \label{zh-11-27-3}
\end{equation}
leads a reduction of $F$ to the form \eqref{zh-7} with $A_{m,N}=0$.

By a direct calculation we obtain
\begin{align}
\frac{\p^3\tilde{F}}{\p\theta^N\p\theta^p\p\theta}=&\sum_{s=1}^{N-p+1} \left(\binom{p+s-1}{s}+\frac12 \binom{p+s-1}{s-1}\right) u^s \frac{\p x_0}{\p u^{p+s-1}}-\frac{\p y_p}{\p u^{N-1}}\nn\\
&-\delta_{p,m+1} \left(N+\frac12\right) u^1\frac{\p x_2}{\p u^{N}},\quad m+1\le p\le N-1,\label{zh-12}\\
\frac{\p^3\tilde{F}}{\p\theta^N\p\theta^m\p\theta}=&\sum_{s=3}^{N-m+1} \left(\binom{m+s-1}{s}+\frac12 \binom{m+s-1}{s-1}\right) u^s \frac{\p x_0}{\p u^{m+s-1}}\nn\\
&+\frac{N-4}2\,\delta_{m,1}\, x_0  +u^1 \left(\frac{2 m+1}2 \frac{\p x_0}{\p u^m}- \frac{2N+1}2 \frac{\p x_1}{\p u^N}+B_1\right)\nn\\
&+u^2 \left(\frac{(m+1)^2}2 \frac{\p x_0}{\p u^{m+1}}+B_2\right)-\frac{\p y_m}{\p u^{N-1}}.\label{zh-10}
\end{align}

We first proceed to fix the differential polynomials  $x_s=x_s(u,u^1,\dots,u^N),\ s=0,1$ and $y_m$. In order to do so, we first choose a solution $x_0$ of the equation 
\[\frac{(m+1)^2}2\frac{\p x_0}{\p u^{m+1}}=-B_2\] 
so that it also satisfies the following additional requirements 
\begin{equation}\label{zh-9}
\frac{\p^2 x_0}{\p u^{m+2}\p u^N}= \dots=\frac{\p^2 x_0}{\p u^{N}\p u^N}=0.
\end{equation}
This can be achieved due to the properties of $B_2$ given in \eqref{zh-8}.
With such a choice of the function $x_0$, the first term that appears in the r.h.s. of \eqref{zh-10} does not depend on $u^N$,
so it can be canceled by the last term $-\frac{\p y_m}{\p u^{N-1}}$
with an appropriately chosen $y_m\in\A^{(N-1)}$. Now we can fix $x_1$ by the equation
\begin{equation}\label{zh-11}
\frac{2N+1}2 \frac{\p x_1}{\p u^N}=\frac{2m+1}2 \frac{\p x_0}{\p u^m}+B_1.
\end{equation}
Thus the above choice of $x_0, x_1$ and $y_m$ enables us to reduce the r.h.s. of \eqref{zh-10} to zero.

In order to reduce $F$ to the form \eqref{zh-7} with $A_{m,N}=0$, 
let us continue to fix the differential polynomials $x_2$ and 
$y_p,\ p=m+1,\dots, N-1$  by the requirement that the r.h.s. of \eqref{zh-12} vanish.
To this end we first fix $x_2$ by the equation
\begin{equation}
\left(N+\frac12\right) \frac{\p x_2}{\p u^N}=\left(m+\frac32\right) \frac{\p x_0}{\p u^{m+1}}.
\end{equation} 
We then choose the functions $y_p\in\A^{(N-1)}$ in order to cancel the rest terms in the r.h.s. of \eqref{zh-12}. This can be done since $x_0$ satisfies the 
condition \eqref{zh-9}.

\noindent\ ii)\ The $m=1$ case. By using \eqref{divide} we can represent $A_{1,N}$ in the form
\[A_{1,N}=u^1 B_1(u,u^1,\dots,u^N)+B_2(u, u^2,\dots,u^N)+B_3(u).\]
Here $B_2$ either equals to zero or, as a polynomial of $u^2,\dots, u^N$, has degree greater than zero, and it also 
satisfies the condition \eqref{zh-8}.
We modify $F$ as we did in \eqref{zh-11-27-1}--\eqref{zh-11-27-3}, then the  r.h.s. of \eqref{zh-10} has the expression
\[\frac{N-4}2 x_0+\sum_{s=1}^N \left(1+\frac{s}2\right) u^s\frac{\p x_0}{\p u^s}+B_2+u^1\left( B_1-\frac{2 N+1}2 \frac{\p x_1}{\p u^N}\right)+B_3-\frac{\p y_1}{\p u^{N-1}}.\]
By using Lemma \ref{lem-54} (when $N>4$), Remark \ref{rmk-55} (when N=4) and the condition \eqref{zh-8}
we can find a solution $x_0=x_0(u,u^1,\dots,u^n)$ satisfying the  equation
\[\frac{N-4}2 x_0+\sum_{s=1}^N \left(1+\frac{s}2\right) u^s\frac{\p x_0}{\p u^s}+B_2=0\]
and the condition \eqref{zh-9}. We then fix the function $x_1$ and $y_1$ by the equations
\[B_1-\frac{2 N+1}2 \frac{\p x_1}{\p u^N}=0,\quad B_3-\frac{\p y_1}{\p u^{N-1}}=0\]
in order to reduce the r.h.s. of  \eqref{zh-10} to zero. 

As we did in the $m\ge 2$ case, we can continue to fix the functions $x_2$ and $y_2,\dots, y_{N-1}$ by
requiring that the r.h.s. of \eqref{zh-12} vanish, so $F$ is reduced to the form \eqref{zh-7} with $A_{1,N}=0$.

Up to this moment, we have proved the existence of $X, Y\in\AA^2$ such that after the modification $F\mapsto F-D_2 X+D_1 Y$
the trivector $F$ can be reduced to the following form
\begin{equation}\label{zh-15}
F=\sum_{1\le p<q\le N-1} A_{p,q}\theta\theta^p\theta^q+B,
\end{equation}
with $A_{p,q}\in\A^{(N)}$, $B\in\AA^{3,(N-1)}$. Now let us consider the vectors $Z_s$ defined in \eqref{zh-6} again. With the 
above form of $F$ we have the following identities:
\[0=\frac{\p Z_s}{\p \theta^p}=\left(N+\frac12\right) u^1\frac{\p^2 A_{p,N-s}}{\p u^N\p u^N},
\quad 1\le p\le N-s-1,\quad 1\le s\le N-2.\]
From these identities we know that $A_{p,q}$'s depend linearly on 
$u^N$. So by using Lemma \ref{lem-55} we can find $X, Y\in\AA^2$ such that 
$F-D_2 X+D_1 Y\in \AA^{3,(N-1)}$.

By induction on $N$,  we  thus proved that $F$ can indeed be reduced to an element of the space $\AA^{3,(3)}$ by modifying it in the way $F\mapsto F-D_2 X+D_1 Y$ with certain $X, Y\in\AA^2$.

\vskip 1em
\noindent Step 2:\ For the case when $N=3$.

In this case, one can also use the above procedure to eliminate the term $\theta\theta^{N-1}\theta^N$ of the trivector $F$,
so we assume that
\[F=A_{1,2}\theta\theta^1\theta^2+A_{1,3}\theta\theta^1\theta^3+B\]
with $A_{1,2}, A_{1,3} \in \A^{(3)}$, $B\in\AA^{3,(2)}$. 
The vanishing of $D_1 D_2 F$ is equivalent to the condition
\begin{equation}
7 u^1\frac{\p^2 A_{12}}{\p u^3\p u^3}-9 u^2 \frac{\p^2 A_{13}}{\p u^3\p u^3}
-5 u^1 \frac{\p^2 A_{13}}{\p u^2\p u^3}=0, \label{d1d2f}
\end{equation}
which implies that
\[u^1\left|\frac{\p^2 A_{13}}{\p u^3\p u^3}\right.,\]
so $A_{1,3}$ can be represented in the following form
\begin{equation}\label{zh-16}
A_{1,3}=u^1 (u^3)^2 b_1+u^3 b_2+b_3,\  b_1\in\A^{(3)},\quad  b_2, b_3\in\A^{(2)}.
\end{equation}
We proceed to reduce $F$ by using the bivector $X, Y\in\AA^2$ of the form
\[X=x_1\theta\theta^1+x_2\theta\theta^2+x_3\theta\theta^3,\quad Y=u X-y_1\theta\theta^1-y_2\theta\theta^2,\]
where $x_1, x_2, x_3\in\A^{(3)}$, $y_1, y_2\in\A^{(2)}$.
The coefficients of $\theta\theta^1\theta^3$ and $\theta\theta^2\theta^3$ in $F-D_2 X+D_1 Y$ are given by
\begin{align}
&R_1=A_{1,3}-\frac72u^1\frac{\p x_1}{\p u^3}-\frac12 x_3+\frac52 u^3 \frac{\p x_3}{\p u^3}+
2 u^2 \frac{\p x_3}{\p u^2}+\frac32 u^1 \frac{\p x_3}{\p u^1}-\frac{\p y_1}{\p u^2},\\
&R_2=\frac92 u^2 \frac{\p x_3}{\p u^3}+\frac52 u^1 \frac{\p x_3}{\p u^2}
-\frac72 u^1 \frac{\p x_2}{\p u^3}-\frac{\p y_2}{\p u^2}
\end{align}
respectively. Let $c\in\A^{(2)}$ be the unique solution of the equation
\[b_2+2 c+2 u^2\frac{\p c}{\p u^2}=0,\]
then one can verify that
\begin{align*}
&x_3=c\,u^3,\ x_2=\frac57\p^{-1}_{u^3}\left(\frac{\p x_3}{\p u^2}\right),\ 
x_1=\frac27\p^{-1}_{u^3}\left(b_1(u^3)^2+\frac32\frac{\p x_3}{\p u^1}\right),\\
&y_2=\frac92\p^{-1}_{u^2}\left(c\,u^2\right),\ y_1=\p^{-1}_{u^2}b_3
\end{align*}
satisfy the equation $R_1=R_2=0$, so the terms in $F$ with $A_{1,3}$ are eliminated.

Now we can use the condition \eqref{d1d2f}
to derive the linear dependence of $A_{1,2}$ on $u^3$, then, by using Lemma \ref{lem-55},
 $F$ is reduced to an element of  the space $\AA^{3,(2)}$. 

\vskip 1em
\noindent Step 3:\ For the case when $N=2$.

In this case $F$ must take the following form
\[F=A\,\theta\theta^1\theta^2,\quad  A\in\A^{(2)}.\]
It is easy to see that this $F$ satisfies $D_1 D_2 F=0$ automatically. We are to find
\[X=x_1\,\theta\theta^1+x_2\,\theta\theta^2,\quad x_1, x_2\in\A^{(2)}\]
such that
\begin{align*}
0=&D_1 F-D_1 D_2 X\\
=&\left(-\frac{\p A}{\p u^2}-\frac{\p x_2}{\p u^2}+\frac52 u^1\frac{\p^2 x_1}{\p u^2\p u^2}-2 u^2\frac{\p^2 x_2}{\p u^2\p u^2}-\frac32 u^1  \frac{\p^2 x_2}{\p u^1\p u^2}\right)\theta\theta^1\theta^2.
\end{align*}
We rewrite $A$ as
\[A=b_0(u,u^2)+u^1 b_1(u,u^1,u^2),\]
and take
\[x_1=\frac25 \p^{-1}_{u^2} b_1, \quad   x_{2}=h(u,u^2),\]
then we obtain the following equation for the function $h$:
\[-\frac{\p b_0}{\p u^2}=\frac{\p h}{\p u^2}+2 u^2 \frac{\p^2 h}{\p u^2 \p u^2},\]
which has a solution due to Lemma \ref{lem-54}.

We thus proved the theorem.
\end{prf}

\subsection{Proof of Theorem \ref{thm-zh1} and Theorem \ref{thm-zh2}}

With the results of Theorem \ref{thm-6} and Theorem \ref{thm-zh7}, the
 proof of Theorem \ref{thm-zh1} follows directly from the long exact sequence \eqref{bhc-long}.
 
To prove Theorem \ref{thm-zh2}, we first check that the infinitesimal deformation part of \eqref{zh-23},  \eqref{zh-24} is
 indeed given by an element of $BH^2_3(\FF)$. Note that the bivectors $(P_1, P_2)$ corresponding to the bihamiltonian structure
 \eqref{zh-19}, \eqref{zh-20} have been given in \eqref{bh-dkdv}. 
 The bivector corresponding to this infinitesimal deformation part reads
\[(0, Q)=\left(0, \ \frac12\int \left(3\,c(u) \theta\theta^3+\frac92\,c'(u) u^1 \theta\theta^2
+\frac32 \left(c''(u) (u^1)^2+c'(u) u^2 \right)\theta\theta^1\right)\right).\]
We denote $Q$'s density by $\alpha$.
Note that this density is first given by Lorenzoni in \cite{Loren}. Let us explain  where it comes from now.
First, we take a representative
\[\omega=-\frac34c(u)\theta\theta^1\theta^2\]
of an element  $[\omega]$ of $BH^3_3(\AA)$. Since $H^3_3(\AA_\lambda, D_\lambda)$ and $BH^3_3(\AA)$ are isomorphic,
there is an element $[\omega_\lambda]$ of $H^3_3(\AA_\lambda, D_\lambda)$ corresponding to this $[\omega]$, which has
a representative
\[\omega_\lambda=-\frac34\left(c(u)+(u-\lambda)c'(u)\right)\theta\theta^1\theta^2.\]
Then by using the connecting homomorphism of the long exact sequence \eqref{long2}, the density $\alpha$ can be fixed by the condition
$\p \omega_\lambda=D_\lambda \alpha$.
Now we know that $Q$ is a representative of an element of $H^2_3(\FF_\lambda, d_\lambda)$,
so it is also a representative of an element of $BH^2_3(\FF)$.

After obtained the infinitesimal deformation $(0, Q)$, we can construct the higher order deformations by using the method
introduced in \cite{LZ1, AL}. More precisely, we want to find a deformation
\[(P_1,\ P_2+\e^2 Q_1+\e^4 Q_2+\cdots+)\]
with $Q_1=Q$. Suppose we have already found $Q_1, \dots, Q_{m-1}$, which satisfy
\begin{equation}
[P_1, Q_k]=0,\quad [P_2, Q_k]=-\frac12\sum_{i=1}^{k-1}[Q_i, Q_{k-i}], \quad k=1, \dots, m-1. \label{ky-13}
\end{equation}
Then we need to find $Q_m$, which satisfies
\begin{equation}
[P_1, Q_m]=0,\quad [P_2, Q_m]=-\frac12\sum_{i=1}^{m-1}[Q_i, Q_{m-i}]. \label{ky-11}
\end{equation}
According to the triviality of $H^2_{2 m+1}(\FF, P_1)$, any solution $Q_m$ of the first equation in \eqref{ky-11}
has the following form
\[Q_m=[P_1, X_m],\quad \mbox{where }X_m\in\FF^1,\]
so we only need to find a solution $X_m$ of the following equation
\begin{equation}
[P_1, [P_2, X_m]] =\frac12\sum_{i=1}^{m-1}[Q_i, Q_{m-i}], \label{ky-12}
\end{equation}
which is equivalent to the second equation of \eqref{ky-11}.

We denote by $W_m$ the right hand side of \eqref{ky-12}, then it is easy to see that
\[[P_1, W_m]=0.\]
By using the equations \eqref{ky-13}, we also have
\begin{align*}
[P_2, W_m]=&\frac12\sum_{k+l=m}[P_2, [Q_k, Q_l]]\\
=&-\frac12\sum_{k+l=m}\left([[P_2, Q_k], Q_l]+[Q_k, [P_2, Q_l]]\right)\\
=&-\sum_{k+l=m}[[P_2, Q_k], Q_l]\\
=&\frac12\sum_{k+l=m}\,\sum_{i+j=k}[[Q_i, Q_j], Q_l]\\
=&\frac16\sum_{i+j+l=m}\left([[Q_i, Q_j], Q_l]+[[Q_j, Q_l], Q_i]+[[Q_l, Q_i], Q_j]\right)\\
=&0,
\end{align*}
so $W_m\in \FF^3_{2 m+2}\cap\Ker d_1 \cap \Ker d_2$. According to Definition \ref{ky-15} of $BH^3_d(\FF)$
and its triviality when $d\ge6$ (see Theorem \ref{thm-zh1}), there exists $X_m\in \FF^1$ such that
\[W_m=[P_1, [P_2, X_m]].\]
The Theorem \ref{thm-zh2} is proved.

\section{Some examples}\label{sec-6}
In this section, we give some concrete examples of bihamiltonian structures 
of the form \eqref{zh-23}, \eqref{zh-24}. We will now denote the $x$-derivatives of the function $u=u(x)$ by $u_x, u_{xx}, u^{(3)}$ and so on,
and, unlike in the previous sections, we denote by $u^m$ the $m$-th power of $u$.

The existence of bihamiltonian structures of the form \eqref{zh-23}, \eqref{zh-24} was
first conjectured by Lorenzoni in \cite{Loren}, and he wrote down its approximation up to $\e^4$,
which are given by the following expressions:
\begin{align}
A_{2,0}=&9 c(u) c'(u),\quad A_{2,1}=\frac{45}2 \left[c'(u)^2+c(u) c''(u)\right] u_x,\nn\\
A_{2,2}=&27 \left[c'(u)^2+c(u) c''(u)\right] u_{xx}+18 \left[3 c'(u) c''(u)+c(u) c^{(3)}(u)\right] u_x^2,\nn\\
A_{2,3}=&18 \left[c'(u)^2+c(u) c''(u)\right] u^{(3)}+27\left[3 c'(u) c''(u)+c(u) c^{(3)}(u)\right] u_x u_{xx}\nn\\
&+\frac92\left[3 c''(u)^2+4 c'(u) c^{(3)}(u)+c(u) c^{(4)}(u)\right] u_x^3,\nn\\
A_{2,4}=&\frac92 \left[ c'(u)^2+c(u) c''(u)\right] u^{(4)}+9 \left[3 c'(u) c''(u)+c(u) c^{(3)}(u)\right] u_x u^{(3)}\nn\\
&+\frac92 \left[3 c'(u) c''(u)+c(u) c^{(3)}(u)\right] u_{xx}^2\nn\\
&+\frac92 \left[3 c''(u)^2+4 c'(u) c^{(3)}(u)+c(u) c^{(4)}(u)\right] u_x^2 u_{xx}.\nn
\end{align}
In a subsequent paper \cite{AL}, Arsie and Lorenzoni extend the approximation up to $\e^8$.

For any given smooth function $g(u)$ define the function $f$ by
\[f(u)=\p_u^{-2} \left(u g''(u)+\frac12 g'(u)\right).\]
Then we have the following bihamiltonian equation
\begin{equation}\label{zh-50}
\frac{\p u}{\p t}=\{u(x), H_f\}_1=\{u(x),H_g\}_2.
\end{equation}
Here the Hamiltonians are given by
\begin{equation}\label{zh-51}
H_f=\int_{S^1}\left[f(u)-\e^2 c(u) f^{(3)}(u) u_x^2+\e^4 \left(A_1 u_{xx}^2+A_2 u_x^4\right)\right] dx+\mathcal{O}(\e^6)
\end{equation}
with
\[A_1=3 c(u) c'(u) f^{(3)}(u)+\frac65 c(u)^2 f^{(4)}(u),\]
\[A_2=-\left[\frac12 c'(u)^2 +c(u) c''(u)\right] f^{(4)}-c(u) c'(u) f^{(5)}(u)-\frac16 c(u)^2 f^{(6)}(u).\]
The Hamiltonian $H_g$ is obtained by replacing $f$ with $g$ in $H_f$. In particular,
if we set 
\[f(u)=\frac{u^{p+2}}{(p+2)!},\quad g(u)=\frac{2 u^{p+1}}{(2p+1) (p+1)!},\quad p=0,1,\dots\]
then \eqref{zh-50} gives the bihamiltonian deformation of the dispersionless KdV hierarchy \eqref{zh-31}.

When $c(u)=\frac1{24}$ the bihamiltonian structure given in Theorem \ref{thm-zh2}
has the following truncated form
\begin{align}
&\{u(x),u(y)\}_1=\delta'(x-y), \label{zh-26}\\
&\{u(x),u(y)\}_2=u(x)\delta'(x-y)+\frac12 u_x(x) \delta(x-y)+\frac{\e^2}8 \delta'''(x-y) \label{zh-27}
\end{align}
which gives the bihamiltonian structure \cite{ZF, Magri-1} for the prototypical nonlinear integrable evolutionary PDE -- the KdV equation
\[\frac{\p u}{\p t}=u u_x+\frac{\e^2}{12} u_{xxx}\]
and the associated KdV hierarchy. 

When $c(u)=\frac{u}{24}$, the bihamiltonian structure \eqref{zh-23}, \eqref{zh-24} is not truncated, however it
is equivalent, under a Miura type transformation $u\mapsto \tilde{u}=u-\frac{\e^2}{16} u_{xx}-\frac{\e^4}{512} u^{(4)}+\mathcal{O}(\e^6)$, to the following bihamiltonian structure:
\begin{align}
&\{\tilde{u}(x),\tilde{u}(y)\}_1=\delta'(x-y)-\frac{\e^2}{8}\delta'''(x-y), \label{zh-28}\\
&\{\tilde{u}(x),\tilde{u}(y)\}_2=\tilde{u}(x)\delta'(x-y)+\frac12 \tilde{u}_x(x) \delta(x-y). \label{zh-29}
\end{align}
it gives the bihamiltonian structure of the equation
\begin{equation}\label{1-ch}
(v-\frac{\e^2}{8} v_{xx})_{t}=v v_x-\frac{\e^2}{12} v_x v_{xx}-\frac{\e^2}{24} v v_{xxx}
\end{equation}
with 
\begin{equation}
\tilde{u}=v-\frac{\e^2} 8 v_{xx}.
\end{equation}
In fact, if we define the Hamiltonians 
\[H_1=\int_{S^1}\left(\frac16\, v^3+\frac{\e^2}{48}\,v v_x^2\right) dx,\quad 
H_2=\int_{S^1}\left(\frac13\,v^2+\frac{\e^2}{24} v_x^2\right) dx,\]
then the equation \eqref{1-ch} can be represented as
\[\tilde{u}_t=\{\tilde{u}(x), H_1\}_1=\{\tilde{u}(x), H_2\}_2.\]
This equation is equivalent to the well-known 
Camassa-Holm equation which, like the KdV equation, describes the dynamics of shallow water waves \cite{FF,  CH, CHH, Fokas, Fu}.

Although the bihamiltonian structures \eqref{zh-23}, \eqref{zh-24} and the associated integrable hierarchies are represented as infinite power series in the dispersion parameters $\e$,  we hope that 
apart from the above two well-known examples of the KdV and Camassa-Holm  integrable hierarchies these new bihamiltonian hierarchies will find applications. 

In the rest part of this section, let us study in some details the bihamiltonian structure  \eqref{zh-23}, \eqref{zh-24} with central invariant $c(u)$ being inversely proportional to $u$. We are to give evidence that it is actually related to the integrable hierarchy obtained from the KdV hierarchy by certain reciprocal transformation.  
To simplify the notations, we will call the family of bihamiltonian equations 
given by \eqref{zh-50}, \eqref{zh-51} with $c(u)=\frac1{24}$ the KdV hierarchy. At the approximation up to $\e^2$ these equations take the form
\begin{equation}
\frac{\p u}{\p t}=f''(u)\,u_x+\frac{\e^2}{24}\left(2 f^{(3)}(u)u^{(3)}
+4 f^{(4)}(u)u_x u_{xx}+f^{(5)}(u)u_x^3\right)+\mathcal{O}(\e^4). \label{zh-52}
\end{equation}
Let us take $f$ to be a specific smooth function $\rho(u)$  (whose explicit definition will be given later) , and denote by $\tau$ the time variable of the 
resulting bihamiltonian evolutionary equation \eqref{zh-50}. Denote $h(u)=\rho''(u)$, then this equation has the expression
\begin{equation}
\frac{\p u}{\p \tau}=h(u)\,u_x+\frac{\e^2}{24}\left(2 h'(u)u^{(3)}+4 h''(u)u_x u_{xx}+h^{(3)}(u)u_x^3\right)+\mathcal{O}(\e^4). \label{zh-53}
\end{equation}
Now we perform the reciprocal transformation to the equation \eqref{zh-52} by
replacing the spatial variable $x$ with the time variable $\tau$ of \eqref{zh-53}.
In order to do so we first obtain from \eqref{zh-53} the following relation:
\begin{align}
u_{x}=&\frac{u_{\tau}}{h(u)}+\frac{\e^2}{24}\left(
-\frac{2 u_{\tau\tau\tau} h'(u)}{h(u)^4}-\frac{4 h''(u) u_{\tau}u_{\tau\tau}}{h(u)^4}
+\frac{18 h'(u)^2 u_{\tau} u_{\tau\tau}}{h(u)^5}\right.\nn\\
&\left.-\frac{h^{(3)}(u) u_{\tau}^3}{h(u)^4}-\frac{24 h'(u)^3 u_{\tau}^3}{h(u)^6}
+\frac{14 h'(u) h''(u) u_{\tau}^3}{h(u)^5}\right)
+\mathcal{O}(\e^4), \label{eq-kdv-xt}
\end{align}
which enables us to represent the equation \eqref{zh-52} in the form
\begin{align}
&\frac{\p u}{\p t}=p(u)\, u_\tau+\frac{\e^2}{12\,h(u)^2}\left[ p'(u)\, u_{\tau\tau\tau}+\left( 2 p''(u)-\frac{5\, h'(u) \,p'(u)}{h(u)}\right) u_\tau u_{\tau\tau} \right.\nn\\
& \quad +\left.\left(\frac{4 \,h'(u)^2\, p'(u)}{h(u)^2}-\frac{5\,h'(u)\, p''(u)}{2\,h(u)}-
\frac{3\,h''(u) \,p'(u)}{2\,h(u)}+\frac12\,p'''(u) \right) u_\tau^3\right]+\mathcal{O}(\e^4),\label{zh-54}
\end{align}
where $p(u)=\frac{f''(u)}{h(u)}$.
In \cite{XZ} it is proved that for a bihamiltonian system of hydrodynamic type, 
after a reciprocal transformation of the above form the transformed system
preserves its bihamiltonian property. So the leading term of the above equation \eqref{zh-54} has a bihamiltonian structure. We conjecture that for bihamiltonian equations with $\e$-deformations a
reciprocal transformation of the above form still preserves the bihamiltonian property. For our particular example of \eqref{zh-54}, we can employ an algorithm proposed in \cite{LZ3, LZ2} to obtain its bihamiltonian structure
at the approximation up to a certain order of $\e$. Here we did the computation at the approximation up to  $\e^8$. Now let us choose the function $h(u)$ carefully, so that the leading terms of the bihamiltonian structure of \eqref{zh-54} coincide, after certain 
change of coordinates, with \eqref{zh-19}, \eqref{zh-20}. To this end, we take
\[ h(u)=u^{-\frac32}.\]
Then after the change of coordinate
\[w=\frac1{u}\]
the equation \eqref{zh-54} has the following bihamiltonian structure
\begin{equation}
\frac{\p w}{\p t}=\{w(\tau), H_1\}_1=\{w(\tau), H_2\}_2,
\end{equation}
with
\begin{align}
&\{w(\tau), w(\sigma)\}_1=\delta'(\tau-\sigma)\nn\\
&\quad 
+\e^2\left[\frac{5}{8\,w(\tau)^2}\delta^{(3)}-\frac{15\,w'(\tau)}{8\,w(\tau)^3}\,\delta''+\left(\frac{3\,w'(\tau)^2}{16\,w(\tau)^4}+\frac{w''(\tau)}{2\,w(\tau)^3}\right) \delta'\right.\nn\\
&\quad \quad +\left.\left(\frac{27\,w'(\tau)^3}{8\,w(\tau)^5}-\frac{27\,w'(\tau) w''(\tau)}{8\, w(\tau)^4}+\frac{9\,w'''(\tau)}{16\,w(\tau)^3}\right)\delta\right]+
\mathcal{O}(\e^4),\label{zh-56}
\\
&\{w(\tau), w(\sigma)\}_2=w(\tau)\delta'(\tau-\sigma)+\frac12\,w'(\tau)
\delta(\tau-\sigma)\nn\\
&\quad +\e^2\left[\frac{1}{2 w(\tau)}\delta^{(3)}-\frac{3\, w'(\tau)}{4\, w(\tau)^2}\,\delta''+\left(-\frac{47\, w'(\tau)^2}{32\, w(\tau)^3}+\frac{7\, w''(\tau)}{8\, w(\tau)^2}\right) \delta'\right.\nn\\
&\quad \quad +\left.\left(\frac{189\, w'(\tau)^3}{64\,w(\tau)^4}-\frac{99\,w'(\tau) w''(\tau)}{32\,w(\tau)^3}+\frac{9\, w'''(\tau)}{16\,w(\tau)^2}\right)\delta\right]+\mathcal{O}(\e^4),\label{zh-57}
\end{align}
and the Hamiltonians 
\begin{align*}
&H_1=\int_{S^1}\left[q(w(\tau))+\frac{\e^2}2\,w'(\tau)^2 \left(\frac{9\,q'(w(\tau))}{16\, w(\tau)^3}+\frac{5\,q''(w(\tau))}{8\, w(\tau)^2}+\frac{q'''(w(\tau))}{12\,w(\tau)}\right)+\mathcal{O}({\e^4}) \right] d\tau\\
&H_2=\int_{S^1}\left[r(w(\tau))+\frac{\e^2}2\,w'(\tau)^2 \left(\frac{9\,r'(w(\tau))}{16\, w(\tau)^3}+\frac{5\,r''(w(\tau))}{8\, w(\tau)^2}+\frac{r'''(w(\tau))}{12\,w(\tau)}\right) +\mathcal{O}({\e^4})\right] d\tau.
\end{align*}
Here $\delta^{(k)}=\delta^{(k)}(\tau-\sigma)$, and the function $q(w), r(w)$ are defined by the relations
\[p(w)=q''(\frac1{w}),\quad q'(w)=-\frac{r(w)}2+ w\, r'(w).\]
The above bihamiltonian structure \eqref{zh-56}, \eqref{zh-57} is equivalent to the one given in \eqref{zh-23},
\eqref{zh-24} with central invariant $c(u)=-\frac1{24\, u}$ under the
Miura type transformation
\[\tilde{w}=w+\frac{\e^2}{32\,w^3}\left(29\, (w')^2-10\,w\, w''\right)+ \mathcal{O}({\e^4}).\]
The above calculation is actually done at the approximation up to $\e^8$, which 
gives us evidence to conjecture that the evolutionary equation \eqref{eq-kdv-xt}
that is obtained from the KdV hierarchy by the above reciprocal transformation is equivalent, under a Miura
type transformation, to the bihamiltonian hierarchy \eqref{zh-54} with central invariant $c(u)=-\frac{1}{24\,u}$.

A full description of the above integrable hierarchy and its bihamiltonian structure is still missing. The approach given
in \cite{KKVV} may be helpful to this problem.

\section{Conclusion}\label{sec-7}

In this paper we proved the existence of deformations of the bihamiltonian 
structure \eqref{zh-19}, \eqref{zh-20}. We hope that apart from the two 
well-known bihamiltonian structures that are associated to the KdV hierarchy and the Camassa-Holm hierarchy, this class of deformed bihamiltonian structures also contains some new bihamiltonian structures which
have important applications in the theory of integrable systems and in
mathematical physics. We will continue to study the properties and applications of these bihamiltonian structures
in subsequent publications. It is also interesting to compute the bihamiltonian
cohomologies $BH_d^p(\FF)$ for $p\ge 4,\ d\ge 6$ which, as we expected, should be trivial.

In the subsequent paper \cite{DLZ2}, we are to consider the existence problem
for a general semisimple bihamiltonian structure of hydrodynamic type by using our formulation of the infinite dimensional bihamiltonian structures and their
cohomologies given in the present paper.

\vskip 5.8 em
\noindent{\bf Acknowledgments.}
\vskip 0.5em
The authors are grateful to Boris Dubrovin for his encouragements and advises. They also thank the referees
for helpful comments and remarks.
This work is partially supported by the NSFC No.\,11071135, 11222108 and No.\,11171176,  and also by the Marie Curie IRSES project RIMMP.

\end{document}